\documentclass[12pt]{article}

\usepackage{setspace}

\usepackage{endnotes}
\usepackage{amsmath}
\usepackage{enumerate}
\usepackage[utf8]{inputenc}
\usepackage {inputenc}
\usepackage{amssymb}
\usepackage{color,soul}
\usepackage{anysize}
\usepackage{tabto}
\usepackage{enumitem}
\usepackage{mathtools}
\usepackage{multirow}
\usepackage{booktabs}
\usepackage{parskip}
\usepackage{hyperref}
\usepackage{longtable}
\usepackage{paralist}
\usepackage{tikz}
\usetikzlibrary{shapes, arrows}
\usepackage{algorithm}
\usepackage{caption}
\usepackage[noend]{algpseudocode}
\usepackage{amsthm}
\setlist{noitemsep}
\let\footnote=\endnote
\let\enotesize=\normalsize
\def\notesname{Endnotes}%
\def\makeenmark{\hbox to1.275em{\theenmark.\enskip\hss}}
\def\enoteformat{\rightskip0pt\leftskip0pt\parindent=1.275em
  \leavevmode\llap{\makeenmark}}
\renewcommand{\arraystretch}{0.5}

\newcommand{\lp}{\bigg(}
\newcommand{\rp}{\bigg)}

\newcommand{\sr}{\mathbb{R}}

\newcommand{\bx}{{\bf x}}
\newcommand{\by}{{\bf y}}
\newcommand{\bz}{{\bf z}}

\newcommand{\xs}[1]{x_{s_#1}}
\newcommand{\ys}[1]{y_{s_#1}}
\newcommand{\zs}[1]{z_{s_#1}}

\newcommand{\cD}{\mathcal{D}}
\newcommand{\cS}{\mathcal{S}}
\newcommand{\ccD}{\bar{\cD}}

\newcommand{\xd}[1]{x_{d_#1}}
\newcommand{\yd}[1]{y_{d_#1}}
\newcommand{\wdz}[1]{w_{d_#1}}

\newcommand{\xb}{\overline{x}}

\newcommand{\Y}{\mathcal{Y}}
\newcommand{\X}{\mathcal{X}}
\renewcommand{\qed}{$\square$}

\newcommand{\bs}{\setminus}

\newcommand{\fun}[1]{\textbf{#1}}
\newcommand{\funi}[1]{\emph{#1}}

\newcommand{\UB}{\mathit{UB}}
\newcommand{\LB}{\mathit{LB}}

\newcommand{\oquo}{``}

\newtheorem{theorem}{Theorem}
\newtheorem{corollary}{Corollary}
\newtheorem{observation}{Observation}
\newtheorem{remark}{Remark}
\newtheorem{definition}{Definition}
\usepackage{fancyhdr}

\makeatletter
\def\BState{\State\hskip-\ALG@thistlm}
\makeatother

\newlength\myindent
\tikzstyle{step} = [rectangle, rounded corners, minimum width=3cm, minimum height=1cm,text centered, draw=black, fill=red!30]
\tikzstyle{arrow} = [thick,->,>=stealth]

\usepackage{natbib}
 \setcitestyle{numbers}
 \bibpunct[, ]{(}{)}{,}{a}{}{,}%
 \def\bibfont{\small}%
\begin{document}



\begin{center}
\textbf{\Large Planar Maximum Coverage Location Problem with Partial Coverage, Continuous Spatial Demand, and Adjustable Quality of Service \\} \vspace{1em}
{ Manish Bansal and Parshin Shojaee\\ Department of Industrial and Systems Engineering, Virginia Tech\\ Email: bansal@vt.edu, parshinshojaee@vt.edu \\ \vspace{1em}
}

\vspace{0.5em}

{First Draft: December 10, 2020}

\end{center}






\noindent \textbf{Abstract.} {%
We consider a generalization of the classical planar maximum coverage location problem (PMCLP) in which \underline{p}artial \underline{c}overage is allowed, facilities have adjustable quality of service (QoS) or service range, and demand zones and service zone of each facility are represented by two-dimensional spatial objects such as rectangles, circles, polygons, etc. We denote this generalization by PMCLP-PC-QoS. A key challenge in this problem is to simultaneously decide position of the facilities on a continuous two-dimensional plane and their QoS. We present a greedy algorithm and a pseudo-greedy algorithm for it, and showcase that the solution value corresponding to the greedy (or pseudo-greedy) solution is within a factor of $1 - 1/e$ (or $1 - 1/e^{\eta}$)  of the optimal solution value where $e$ is the base of natural logarithm and $\eta \leq 1$. We also investigate theoretical properties and propose exact algorithms for solving: (1) 
PMCLP-PC-QoS where demand and service zones are represented by axis-parallel rectangles (denoted by PMCLP-PCR-QoS), which also has applications in camera surveillance and satellite imaging; and (2) one dimensional PMCLP-PC-QoS which has applications in river cleanups.
These results extend and strengthen the only known exact algorithm for PMCLP-PCR-QoS with fixed and same QoS by Bansal and Kianfar [INFORMS Journal on Computing 29(1), 152-169, 2017]. We present results of our computational experiments conducted to evaluate the performance of our proposed exact and approximation algorithms.

%
}



%


\section{Introduction}

Over the years, many classes of problems related to locating service facilities to cover the demand for service have been studied; refer to \cite{ChuMur13,DreHam04} for extensive reviews. 
A well-known classical facility location problem is Maximum Coverage Location Problem (MCLP) for locating $p$ service facilities having known service range with the objective of maximizing the total covered demand \citep{ChuRev74}. In literature, various forms of MCLP consider a finite set of pre-determined candidate positions for the facilities \citep{ChuRev74, MurOKe02, Daskin89}. However, a generalization of the classical MCLP, referred to as the planar MCLP (PMCLP), considers locating the facilities anywhere in a continuous two-dimensional plane ~\citep{Chu84, WatGan82, Mehrez82}. These problems have received considerable attention in the literature due to their widespread applicability ranging from locating emergency healthcare centers, locating fire fighting stations, and making policy through geographical informative systems (GIS) to solve clustering problems which themselves find applications in data mining, machine learning, and bio-informatics \citep{Chu86,FarAsgHeiHosGoh12,SchJayBar93}. 
With the growing availability of spatial data in the foregoing domains, the state-of-the-art facility location analytical tool set needs to evolve. In this paper, we consider a generalization of the PMCLP, by utilizing spatial representation of demand and service zones of facilities, allowing adjustable quality of service (QoS) or service range for the facilities, and allowing partial coverage in its true sense.

In most of the PMCLP literature, demands are represented as aggregated points, which are obtained by aggregating the demand of each zone at a single representative point (e.g., its centroid). \cite{MurTon07} considered PMCLP where the demand zones are represented as line segments or polygons instead of limiting them to be represented by aggregated points. However, similar to the PMCLP with point representation of demand, they assume that demand zones (points, line segments, or polygons) are either completely covered or not covered at all by any service zone. This coverage assumption is referred to as ``binary coverage" \citep{BanKia17IJOC}. The binary coverage is a simplifying assumption to model PMCLP as binary programs, but it 
ignores partial coverage of demand zones. 
%
Researchers have studied its impact and referred to it as region misrepresentation coverage error (RMCE) \citep{CurSch90,TonChu12,Mur16}. These analysis represent that the modeling solutions are sensitive to how demand is represented in the PMCLP \citep{MurOKe02}, and there are evidences that the point representation introduces unintended measurement and interpretation errors \citep{CurSch90,TonMur09} which lead to gaps in actual coverage. 
 
\begin{figure}[!h]
\centering
\includegraphics[width=1\textwidth]{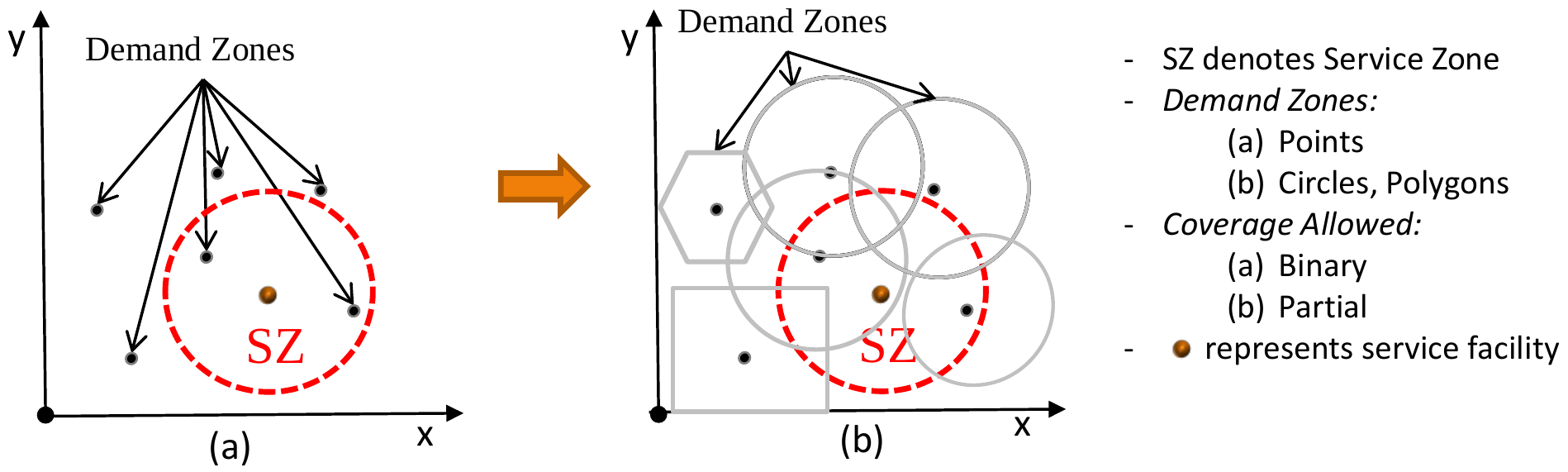} 
\caption{(a) Classical Planar MCLP; (b) Planar MCLP with partial coverage and general representation of demand and service zones. }
\label{fig:RotatedDZs}
\end{figure}

In Figure \ref{fig:RotatedDZs}(a), we present an example of PMCLP where demand zones are represented using points and the service zone of a facility with fixed  service range $r$, using Euclidean distance, is a circle of radius $r$ centered at the facility. Observe that because of the binary coverage assumption, only two demand zones are completely covered by the service zone, whereas others are not at all covered. This causes RMCE because in reality, as shown in \ref{fig:RotatedDZs}(b), five demand zones (represented by spatial objects) are partially covered by the service zone. Note that the service zone of a facility, using the rectilinear distance, is a diamond (a square rotated $45$ degrees) with a diagonal of $2d$ centered at the facility. Lately, \cite{BanKia17IJOC} developed the first exact algorithm for the PMCLP where ``partial coverage" is allowed, and service and demand zones are defined by axis-parallel rectangles. It generalizes PMCLP with rectilinear distance. However, they assumed that all facilities have fixed and same service range or QoS, i.e., dimensions of the rectangular service zones are fixed and same. 

Although attempts have been made to separately address non-point representation of demand, partial coverage, and adjustable QoS, but to our knowledge, no study has tackled them together. Few studies have considered MCLP with gradual coverage, i.e., the coverage level of a demand zone depends on their distance from the facilities \citep{BerKra02,ChuRob83,BerKraDre03}, but similar to MCLP, the demand zones are still represented by points and there is a finite set of pre-specified candidate positions for the facilities. This problem is referred to as gradual coverage location problem. \cite{DreWesDre04} considered planar version of the forgoing problem with only single facility. (See section~\ref{sec:challenges} for more details.) However so far, no study has considered PMCLP with ($i$) general representation of demand and service zones using two-dimensional spatial objects, ($ii$) facilities having adjustable QoS or service range, and ($iii$) partial coverage in its true sense, i.e., when covering only part of a demand zone is allowed and the coverage accrued in the objective function as a result of this is proportional to the demand of the covered area only. We denote this generalization of PMCLP by PMCLP-PC-QoS.

\subsection{Problem definition: PMCLP-PC-QoS}\vspace{-0em}

We define the PMCLP-PC-QoS as follows. Let $\mathcal{D} = \{d_i, i =1,\ldots, n \}$ be a set of $n$ (possibly overlapping) spatial objects, referred to as demand zones (DZs), on a two-dimensional plane such that dimensions and location of each DZ are known. Now, consider another set of spatial objects that provide coverage of facilities, referred to as service zone (SZ), and denote the set of SZs by $\mathcal{S}:= \{s_j, j = 1,\ldots,p\}$. We assume that shape and orientation of all SZs are same and known, but their location and dimensions are unknown. Let $s_0$ be a spatial object with same shape and orientation as of the SZs, but with known dimensions. We refer to $s_0$ as \oquo base" SZ because it provides base reward rate (per unit area) $v_i \in \sr_+$ for covering each DZ $d_i$. Dimensions of each SZ $s \in \mathcal{S}$ is a scalar $z_s \geq 1$ multiple of the dimensions of $s_0$. We refer to $z_s$ as the scaling factor of SZ $s$. The base reward rate is utilized to compute the reward rate $v_i^z = v_i/\eta(z) \leq v_i$ for each DZ $d_i$ covered by a SZ with scaling factor $z$, where $\eta(z)$ is a strictly increasing function. When scaling factor $z$ increases, the dimensions of a SZ increases but the corresponding reward rate decreases. Therefore, by increasing the scaling factor $z$, the service range of a facility increases, but the QoS decreases. We assume that $\Xi_s:=\{\xi^s_1, \ldots, \xi^s_m\}$ is a finite set of $m$ possible scaling factors for SZ $s\in \mathcal{S}$, which allows adjustable QoS for the SZs. In this paper, we will use terms QoS, scaling factor, or service range interchangeably for a SZ $s$ and associate all of them with $z_s \in \Xi_s$ for $s \in \mathcal{S}$. 

The goal of PMCLP-PC-QoS is to identify the location and QoS of all SZs, i.e., $(x_{s_j}, y_{s_j}, z_{s_j}) \in \sr^2 \times \Xi_s$ for all $j = 1,\ldots, p$, where $(x_{s_j}, y_{s_j})$ denotes coordinates of either center or a corner (if exists) of SZ $s_j$. The objective function of this problem is to maximize the total reward captured by these $p$ SZs, i.e., 
\vspace{-0.5em}%
\begin{align} \label{eq:MaxRewardFunc}
\max_{\bx,\by,\bz}\bigg\{ f(\bx,\by,\bz):=  \sum_{i=1}^n f_i(\bx,\by,\bz) = \mathcal{T}_i\lp d_i \cap \lp \cup_{j=1}^p s_j \rp\rp \bigg\}
\end{align} 
where $\cup$ and $\cap$ denote the union and intersection of coverage zones, $\mathcal{T}_i(.)$ returns the total reward captured from DZ $d_i$ by SZs $s_1, \ldots, s_p$, and $(\bx, \by, \bz) = (\xs{1}, \ldots, \xs{p},\ys{1}, \ldots, \ys{p},\zs{1},\linebreak \ldots, \zs{p} )$. 
In case each SZ offers same QoS, i.e., $z_{s_j} = \hat{z}$ for all $j \in \{1,\ldots,p\}$, the function $f(\bx,\by,\bz)$ reduces to $\displaystyle\sum_{i=1}^n v^{\hat{z}}_i \times A\lp d_i \cap \lp \cup_{j=1}^p s_j \rp\rp $ where $A(.)$ returns the area of its argument. 
In this paper, we present greedy and pseudo-greedy approximation algorithms for the PMCLP-PC-QoS 
as well as exact algorithms for its following variants:  
%
%

\noindent (a) 
\emph{PMCLP-PC-QoS with axis-parallel \underline{r}ectangular service and demand zones} (denoted by PMCLP-PC\underline{R}-QoS). Given a set of rectangular DZs identified by coordinates of their lower left corner and their dimensions (width and length), denoted by $(\xd{i},\yd{i})\in \sr^2$ and $(\wdz{i},l_{d_i}) \in \sr^2_+$, respectively, for $i=1,\ldots,n$. We assume that the base SZ $s_0$ is an axis-parallel rectangle with width $w_{s_0}$ and length $l_{s_0}$. Consequently, the SZs $s_1, \ldots, s_p$ are also axis-parallel rectangles. We also assume that the set of scaling factors, $\Xi_s = \Xi:= \{\xi_1, \xi_2, \ldots, \xi_m\}$ for all $s \in \mathcal{S}$. 
The goal of PMCLP-PCR-QoS is to find locations ($\bx,\by$) as well as scaling factor $\bz$ of all SZ such that the maximum reward is covered. Observe that the width and length of SZ $s$ depend on the scaling factor $z_s \in \Xi$ associated to it, i.e., $({w_s^z},{l_s^z}) = (z_s w_{s_0}, z_s l_{s_0})$. A motivation behind studying this problem is its application in telerobotics, camera surveillance, and satellite imaging (refer to the next section for details). Also, the problems studied by \cite{SonStaGol06} and \cite{BanKia17IJOC} are special cases of PMCLP-PCR-QoS where $p=1$ and $\Xi=\{1\}$, i.e., the facilities cannot adjust QoS and their SZs have fixed and same dimensions, respectively. \cite{BanKia17IJOC} also proved that if $p$ is a part of input, PMCLP-PCR-QoS with $\Xi = \{1\}$, denoted by PMCLP-PCR, is NP-hard. Therefore, PMCLP-PC-QoS and PMCLP-PCR-QoS are  also NP-hard because PMCLP-PCR is a special case of these problems.

\noindent (b) \emph{One-dimensional PMCLP-PC-QoS} (1D-PMCLP-PC-QoS): As the name suggests, the one-dimensional facility location problem has DZs and SZs placed on a line. Specifically, we assume that the DZs are line segments on $x$-axis whose coordinate of left corner or end point, $\xd{i}$, $i=1,\ldots,n$, and width, $w_{d_i}$, $i=1,\ldots,n$, are known. Likewise, the SZs are also defined as line segments on $x$-axis with different QoS, i.e., $\Xi_s = \{\xi^s_1\}$. The motivation behind studying this special case is its application in locating trash booms (facilities) for cleaning trash zones (or DZs) in a river (represented by 1D line). In the next section, we discuss more about this application as well. Another application of this problem is in locating $p$ regional wastewater treatment plants on a river that has a given set of fixed segments (DZs) with high priority (reward rate). 
To our knowledge, the 1D-PMCLP-PC-QoS has not been studied before. In literature, one-dimensional plant location problem is studied by \cite{BR1998}, and \cite{Brimberg2001} where DZs are represented by points, and facilities are uncapacitated and capacitated, respectively. As a result, the authors were able to provide mixed binary programs and dynamic programs, respectively, for them. They also show that under certain conditions, these problems are polynomially solvable using linear programs.

\begin{figure}[!h]
\centering
\includegraphics[width=\textwidth]{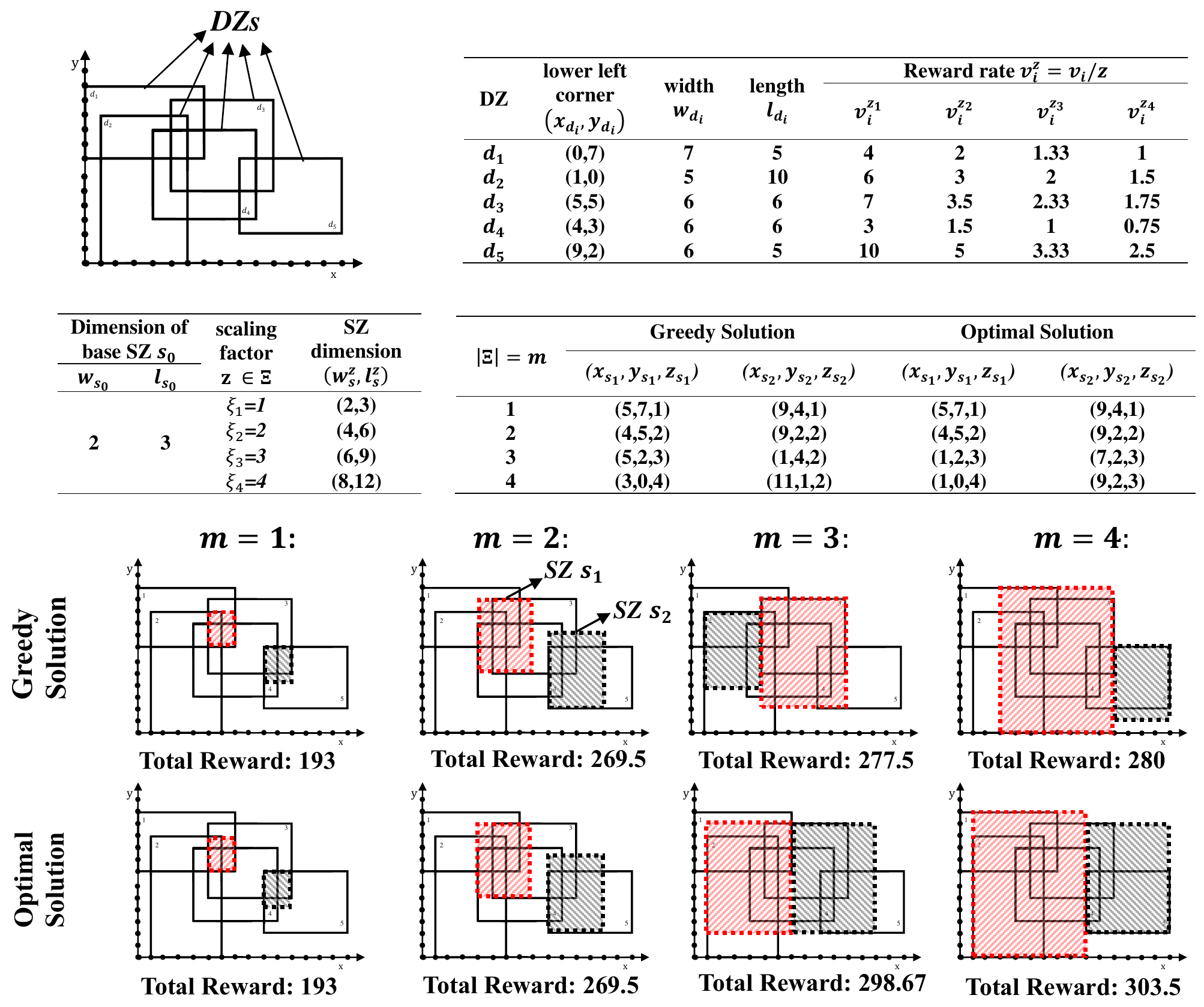} 
\caption{An example of PMCLP-PC-QoS with a set of five axis-parallel rectangular DZs $\{d_1, d_2, \ldots, d_5\}$ with known location, width, length, and base reward rates, a base SZ $s_0$, and a set of scaling factors $\Xi = \{\xi_1, \ldots, \xi_m \}$ for two SZs, i.e. $p=2$. The optimal and greedy solutions for $m \in \{1, \ldots, 4\}$ are obtained from algorithms presented in this paper.
}
\label{fig10}
\end{figure}

 In Figure \ref{fig10}, we provide an example to illustrate PMCLP-PCR-QoS with five DZs $\{d_1, d_2, \linebreak\ldots, d_5\}$ with known location, width, length, and base reward rates, a base SZ $s_0$ with known dimensions, and a set of scaling factors $\Xi = \{\xi_1, \ldots, \xi_m \}$ for two SZs, i.e. $p=2$. A similar example for $m=1$ has been considered by \cite{BanKia17IJOC}. The locations and dimensions are considered integers just for the simplicity of representation. In this example, we assume that $\eta(z) = z$ and therefore, $v^z_i = v_i/z$. Using the algorithms presented in this paper, we present optimal and greedy solutions (along with solution values, i.e., reward) for $m \in \{1, \ldots, 4\}$ and $\xi_k = k$, $k=1,\ldots,4$. Notice that for $m \in \{1,2\}$, the greedy solution is same as the optimal solution for this example. However, for $m \in \{3,4\}$, the optimal solution value is better than the greedy solution value. Moreover, with the increase in $m$, the captured reward is non-decreasing because the solution space for smaller value of $m$ is a subset of solution space for larger value of $m$. This demonstrates a significance of considering facilities with adjustable QoS. 
 
 \subsection{Other Applications of the PMCLP-PC-QoS} \vspace{-0em}

The PMCLP-PC-QoS also has direct applications in variety of emerging domains:

(a) \emph{Telerobotics.} The advent of network-based telerobotic camera systems enable multiple participants or researchers in space exploration, health-care, and distance learning, to interact with a remote physical environment using shared resources. This system of $p$ networked robotic cameras with discrete resolutions receives rectangular requests or subregions (DZs) from multiple users for monitoring \citep{SonStaGol06,XuSonYi10,Xu08}. Each request has an associated reward rate (per-unit area) that may depend on the priority of the user or the importance level associated with monitoring that subregion and the resolution level (QoS) of the camera utilized to cover the subregion. The goal is to select the best view frame for the cameras (rectangular SZs) and their resolution (QoS) to maximize the total reward from the captured parts of the requested subregions. Interestingly, this problem is same as the PMCLP-PCR-QoS. In the literature, \cite{SonStaGol06} studied the PMCLP-PCR-QoS for single camera, i.e. $p=1$, and \cite{XuSonYi10,Xu08} considered the PMCLP with rectangular DZs and binary coverage where the rectangular SZs are not allowed to overlap. The PMCLP-PCR-QoS subsumes the foregoing problems.\vspace{-0.3em}


(b) \emph{River Cleaning.} The main motivation behind studying the 1D-PMCLP-PC-QoS is its application in locating trash-booms for river cleanup. Rivers are the primary recipient of stormwater and as a result, the amount of trash, floating debris, and litter in them is rapidly increasing. 
The major sources along the river that contribute to the trash problem include industrial and recreational parks, and tributaries. 
There are numerous on-going efforts to clean rivers in the US and also across the globe, for example, Nile in Egypt, Yangtze in China, Ganga in India, etc. One way to reduce/cover trash in rivers is by appropriately locating trash-booms in and along the rivers. 
Since
the point representation of a trash zone does not take into account the displacement of the trash from source locations. Therefore, an analogous of 1D-PMCLP-PC-QoS is a trash-boom location problem, where a river is represented by a one-dimensional axis, and trash zones (or DZs) and coverage of trash-booms (SZs) are represented by line segments. In this problem, the service range of each trash boom is known but can be different, and the partial coverage is allowed.

\subsection{Organization of this Paper} 

In Section \ref{sec:challenges}, we discuss challenges in solving PMCLP-PC-QoS using well-known approaches for solving (planar) MCLP and its variants. We present a greedy algorithm and a pseudo-greedy algorithm for solving the PMCLP-PC-QoS (Section~\ref{sec:Greedy}), and showcase that the solution value corresponding to the greedy solution is within a factor of $1 - 1/e$ of the optimal solution value where $e$ is the base of natural logarithm. 
This extends the similar results of \cite{CorFisNem77} and \cite{HocPat98} for special cases of PMCLP-PC-QoS (see Section \ref{sec:LitRevAA} for details). 
In Section \ref{sec:PMCLP-PCR-QoS}, we analyze the objective function and solution space of the PMCLP-PCR-QoS. We strengthen the theoretical properties provided for PMCLP-PCR in \cite{BanKia17IJOC} by reducing the solution search space and thereby, improve the computational efficiency of their algorithm. We also introduce theoretical properties to reduce search space for optimal solution of PMCLP-PCR-QoS, and utilize these properties to develop a branch-and-bound based exact algorithm for it. In Section~\ref{sec:1D-PMCLP-PC}, we introduce 1D-PMCLP-PC-QoS and provide theoretical properties for its solution space along with an exact algorithm for it. 
We also conduct computational experiments to evaluate the performance of our exact algorithms and greedy approach for PMCLP-PCR, PMCLP-PCR-QoS, and 1D-PMCLP-PC-QoS 
(Section~\ref{sec:ComExp}). Finally, we provide our concluding remarks 
in Section~\ref{sec:Conclusion}.


\section{Challenges in solving PMCLP-PC-QoS} \label{sec:challenges} \vspace{-0em}

In this section, we present two well-known approaches for solving (planar) MCLP and its variants, and discuss how they cannot be directly utilized for solving the PMCLP-PC-QoS. 

\subsection{Linear Binary Programming} \label{sec:MIP}

The motivation for the binary coverage assumption in the (planar) MCLP (where demand is represented by points, line segments, or polygons) is to make the problem manageable by readily formulating them as a linear binary program (LBP). This is easy to do for the MCLP because of the discrete nature of candidate locations for service facilities (as per the definition). Moreover, even in studies considering a planar setting, i.e., allowing the facilities to be located anywhere in the continuous plane, coverage is still assumed to be binary \citep{MurTon07} because this helps to show that a finite number of potential facility locations, called the circle intersection point set (CIPS) \citep{Chu84} and polygon intersection  point set (PIPS) \citep{MurTon07}, exist which contain an optimal solution to this problem. Thereby resulting in the following well-known LBP for the (planar) MCLP:
\begin{equation}
\max \bigg\{ \sum\nolimits_{i} v_i x_i: \sum\nolimits_{j} a_{ij} y_j \ge x_i \text{ for all } i, \ \
\sum\nolimits_{j} y_j = p, \ \ x_i, y_j \in \{0,1\} \text{ for all } i,j\bigg\},
\end{equation}
where binary variable $x_i=1$ if demand zone $i$ is covered and $x_i=0$ otherwise, variable $y_j=1$ if a service facility is sited at the candidate/PIP/CIP point $j$ and $y_i=0$ otherwise, $v_i$ is the given total demand of demand zone $i$, and $a_{ij}$ is the given binary  value which is 1 if demand zone $i$ is covered by locating a facility at candidate/CIP/PIP point $j$, i.e., distance between point $i$ is no greater than known service range of a facility (denoted by $r$) located at point $j$ or $dist(i,j) \leq r$.

Furthermore in literature, LBP formulations have also been used to tackle so-called gradual coverage location problem (GCLP). Similar to the MCLP, in GCLP, the demand zones are still represented by points but the coverage level depends on their distance from the facilities. So far LBP formulations are known only for the GCLP defined over a finite set of pre-specified candidate positions for the facilities \citep{BerKra02,BerKraDre03,ChuRob83}, and planar GCLP with single facility, i.e. $p=1$ \citep{DreWesDre04}. More specifically, given a set of $n$ demand points and a set $\mathcal{Y}$ of finite number of positions where $p$ facilities can be placed, the GCLP can be stated as follows: $\max_{F \subset \mathcal{Y}}\big\{ \sum_{i=1}^n v_i g_i(\Delta_i(F)): |F| =p\big\}$ where $F$ is a set of locations where facilities are located, $\Delta_i(F) = \min_{j \in F}\{ dist(i,j)\}$ is the minimum distance between demand point $i$ and any facility location in $F$, and $g_i(.) \in [0,1]$ is a pre-defined coverage function. The following LBP formulations have been derived for the GCLP with linear decay coverage function \citep{BerKraDre03,DreWesDre04}, i.e., $g_i(\Delta) = 1- \beta \Delta$ where $\beta >0$ is a constant, or step-coverage function  \citep{BerKra02,ChuRob83}, i.e., $g_i(\Delta) = \delta_k$ if $\Delta \in (r^i_{k-1}, r^i_k]$ for $k =1,\ldots,K$ where $\delta_1 = 1 > \delta_2 >\ldots > \delta_K =0$ (coverage levels) and $r^i_{0} = 0 < r^i_1 = r < r^i_2 < \ldots < r^i_{K}$ (coverage radii):
\begin{align}\label{eq:GCLP}
\max \bigg\{ \sum_{i=1}^n \sum_{j \in \mathcal{Y}} c_{ij} x_{ij} \big| & \ \ a_{ij} y_j \ge x_{ij}, \  \forall i, j; \  \sum_{j \in \mathcal{Y}} a_{ij} x_{ij} \leq 1 \ \forall i; \nonumber \\ &
\sum\nolimits_{j\in \mathcal{Y}} y_j = p, x_{ij}, y_j \in \{0,1\}, \  \forall i,j\bigg\}
\end{align}
where $c_{ij} = v_i g_i(dist(i,j))$, $x_{ij} = 1$ if facility located at point $j$ covers the demand point $i$ and $x_{ij}=0$ otherwise, and  $a_{ij}$ is the given binary value which is 1 if demand at point $i$ can be partially/completely covered by locating a facility at candidate point $j$, i.e., $dist(i,j) < r^i_K$. Note that if $g_i(\Delta) = 1$ for $\Delta \leq r$ and $g_i(\Delta) =0$ for $\Delta >r$, then the GCLP reduces to the MCLP.

PMCLP-PC-QoS with its features of partial coverage, adjustable QoS, and general spatial representation of DZs and SZs is significantly harder to solve compared to the (planar) MCLP problem with binary coverage, even when a demand zone is represented by a line segment or polygon, and the GCLP. Using LBPs to solve the PMCLP-PC-QoS is no more feasible. 

\subsection{Greedy-Based Algorithms}\label{sec:LitRevAA} \vspace{-0em}

 Since even MCLP is an NP-hard problem, approximation algorithms have also been developed in the literature for solving (planar) MCLP. Among them, greedy approximation algorithm is a well-known approach because it requires solving single facility problem for multiple ($p$) times. In this direction, \cite{CorFisNem77} and \cite{HocPat98} provided greedy algorithms for solving variants of MCLP along with their approximation ratios. More specifically, in a seminal paper on locating bank accounts (or facilities) in at most $p$ out of $m$ known cities to cover $n$ clients (or DZs represented by points), \cite{CorFisNem77} considered a variant of MCLP (or GCLP) in which fixed cost of locating accounts is also deducted in the objective function and it is assumed that the coverage function $g_i(dist(i,j)) = \phi_{ij}$ is constant for each pair of client $i$ and facility location $j$. They derived an LBP formulation which is same as \eqref{eq:GCLP} when the fixed costs are zero, and also presented a greedy approach which provides a solution whose value is within a factor of $1 - 1/e$ of the optimal solution value. Clearly with zero fixed costs, this problem is a special case of PMCLP-PC-QoS. 

\cite{HocPat98} extended the results in \cite{CorFisNem77} by considering a so-called maximum $p$-coverage problem (MCP) which is defined as follows: Given a universal set of elements $U$ where each element $i \in U$ has weight $v_i$ associated to it and a class $\mathcal{V}$ of subsets of $U$, the goal is to select $p$ members (or subsets of $U$) from the class $\mathcal{V}$ such that 
the sum of the weights of the elements in the union of these subsets is maximum. We can build its correspondence with facility location problem, in particular PMCLP-PC-QoS, by considering each element of the set $U$ as a DZ and each member of the class $\mathcal{V}$ as a set of DZs that are completely covered by a SZ located at a pre-specified candidate position and have fixed QoS. Observe that the foregoing problem does not allow partial coverage of elements/DZs and restricts selected members to belong to the class $\mathcal{V}$, which is equivalent to SZs to be located at pre-specified candidate positions. This implies that the MCP is a special case of PMCLP-PC-QoS. \cite{HocPat98} provided greedy and pseudo-greedy algorithms for MCP, and showcase that the solution value corresponding to the greedy (or pseudo-greedy) solution is within a factor of $1 - 1/e$ (or $1 - 1/e^{\eta}$)  of the optimal solution value where $e$ is the base of natural logarithm and $\eta \leq 1$. In this paper, we further extend their algorithmic and theoretical results for PMCLP-PC-QoS. Note that because of the binary coverage assumption of MCP and discrete nature of the bank account location problem, the results in \cite{HocPat98} and \cite{CorFisNem77}, respectively, cannot be directly applied on the PMCLP-PC-QoS. 

\section{Approximation Algorithms for PMCLP-PC-QoS} \label{sec:Greedy} 

In this section, we present greedy and pseudo-greedy algorithms for PMCLP-PC-QoS with $\Xi_s = \Xi = \{\xi_1, \ldots, \xi_m\}$ for all $s \in \mathcal{S}$.

\noindent {\bf Greedy Algorithm.} Assuming that there exists an exact algorithm for solving the PMCLP-PC-QoS with $p=1$, referred to as Single SZ Problem (SSP), we solve multiple SSPs in our greedy-based algorithm for the PMCLP-PC-QoS. The pseudocode is presented in Algorithm 1, where for a set of DZs $\mathcal{D}^g_j$, the function \texttt{SingleSZProblem{$(\mathcal{D}^g_j, s_j)$}}, $j \in \{1,\ldots,p \}$, returns the maximum {reward}, $\psi_g^j$, covered by the SZ $s_j$ along with its optimal position $(x^g_{s_j}, y^g_{s_j})$ {and QoS $z^g_{s_j}$ } (Line 4). We initialize the algorithm in Line 2 by setting $\mathcal{D}^g_1=\mathcal{D}$ (the original set of all given DZs). We also use the function \texttt{TrimOut$(d, s_j, x^g_{s_j}, y^g_{s_j}, z^g_{s_j})$} to eliminate the parts of DZ $d$ that are covered by SZ $s_j$ positioned at $(x^g_{s_j}, y^g_{s_j})$ with dimensions $z^g_{s_j}$ times the dimensions of base SZ $s_0$. In Lines 6-8, we create set $\mathcal{D}^g_{j+1}$ for the next iteration by replacing each DZ $d$ in the set $\mathcal{D}^g_j$ with trimmed DZs. We denote the set of trimmed DZs, that replaces DZ $d$, by $T_d$. The summation of the maximum covered reward by calling \texttt{SingleSZProblem$(\mathcal{D}^g_j, s_j)$} over $j \in \{1,\ldots,p \}$ gives a feasible solution and a lower bound on the optimal objective value of the PMCLP-PC-QoS. Algorithm 1 generalizes the greedy-based polynomial-time heuristic of \cite{BanKia17IJOC} for solving the PMCLP-PCR in $O(n^2p^3 m)$, and Theorem \ref{thm:GreedyAlgo} provides an approximation ratio for this special case as well. 

\begin{algorithm}[H] 
  \fontsize{10}{13}\selectfont
\caption{Greedy Algorithm for PMCLP-PC-QoS}
\begin{algorithmic}[1]
\Function{\fun{GreedyAlgorithm}}{{$\mathcal{D}$, $\mathcal{S}$}}
	\State $\mathcal{D}^g_1 := \mathcal{D}; \psi_g := 0$;
	\For{$j =1,\ldots,p$ }
			 \State $(\psi_g^j, x^g_{s_j}, y^g_{s_j}, z^g_{s_j}):=$\texttt{SingleSZProblem{$\lp \mathcal{D}^g_j, s_j \rp$}};  
			 \State $\psi_g \leftarrow {\psi_g} + \psi_g^j$;
	     \For{$d \in \mathcal{D}^g_j$}
			\State{$\mathcal{D}^g_{j+1}\leftarrow\{\mathcal{D}^g_j \setminus d\} \cup$ \texttt{TrimOut$\lp d, s_j, x^g_{s_j}, y^g_{s_j}, {z^g_{s_j}}\rp$};} 
					\EndFor
	   \EndFor
	  \State \Return{$\lp \psi_g,x^g_{s_1},\ldots, x^g_{s_p}, y^g_{s_1},\ldots, y^g_{s_p}, z^g_{s_1},\ldots, z^g_{s_p} \rp$}
\EndFunction
	\end{algorithmic}
\end{algorithm}

\begin{observation} \label{observation1}
In iteration $j \in \{1,\ldots,p\}$ of the greedy algorithm (Algorithm 1), we exactly solve a SSP for $\mathcal{D}^g_j$ set of DZs. Observe that summing the demand covered in the iteration, i.e., $\psi_g^j$, for $p$ times provides an upper bound on the optimal solution value of the PMCLP-PC-QoS with $\Xi_{s_j} = \Xi$, $s_j \in \mathcal{S}$, 
for $\mathcal{D}^g_j$ set of DZs. This is because the summation does not consider overlapping of the SZs.
\end{observation}

\noindent {\bf Pseudo-Greedy Algorithm.} In the greedy algorithm, we assume that an exact algorithm for solving the SSP is known. However, in case this assumption fails and only an $\eta$-approximate algorithm for solving the SSP is known, then we utilize a pseudo-greedy algorithm for solving the PMCLP-PC-QoS. The pseudo-greedy algorithm is same as the greedy algorithm (Algorithm 1) except that the function \texttt{SingleSZProblem} is replaced by function \texttt{$\eta$-ApproxSSP} which returns an approximate demand, $\kappa_r^j$, covered by the SZ $s_j$ that is within a known factor of $\eta$ ($\leq 1$) of the optimal solution value returned by \texttt{SingleSZProblem}. The pseudocode is given in Algorithm 2. Note that for $\eta =1$, the pseudo-greedy algorithm is exactly same as the greedy algorithm. Furthermore, for $\eta <1$, the sets of DZs $\mathcal{D}^g_j$ and $\mathcal{D}^r_j$, $j = 2,\ldots, p$, in the greedy algorithm (Algorithm 1) and pseudo-greedy algorithm (Algorithm~2), respectively, are different.   

\begin{algorithm}[!h]
  \fontsize{10}{13}\selectfont
\caption{Pseudo-Greedy Method}
\begin{algorithmic}[1]
\Function{\fun{PseudoGreedyAlgorithm}}{{$\mathcal{D}$, $\mathcal{S}$}}
	\State $\mathcal{D}^r_1 := \mathcal{D}; \kappa_r := 0$;
	\For{$j =1,\ldots,p$ }
			 \State $(\kappa_r^j, x^r_{s_j}, y^r_{s_j}, z^r_{s_j}):=$ \texttt{$\eta$-ApproxSSP{$\lp \mathcal{D}^r_j, s_j \rp$}};  
			 \State $\kappa_s \leftarrow {\kappa_r} + \kappa_r^j$;
	     \For{$d \in \mathcal{D}^r_j$}
					\State{$\mathcal{D}^r_{j+1} \leftarrow \{\mathcal{D}^r_j \setminus d \} \cup$ \texttt{TrimOut$\lp d, s_j, x^r_{s_j}, y^r_{s_j}, z^r_{s_j}\rp$};} 
					\EndFor
	   \EndFor
	  \State \Return{$\lp \kappa_r, x^r_{s_1},\ldots, x^r_{s_p}, y^r_{s_1},\ldots, y^r_{s_p}, z^r_{s_1},\ldots, z^r_{s_p} \rp$}
\EndFunction
	\end{algorithmic}
\end{algorithm}

\begin{remark}
The greedy (or pseudo-greedy) algorithm depends on \texttt{SingleSZProblem} (or \texttt{$\eta$-ApproxSSP}) and \texttt{TrimOut}, and the number of trimmed DZs generated after each iteration. Till date these functions are not known, except for PMCLP-PCR-QoS with $p=1$ \citep{BanKia17IJOC,SonStaGol06}. However, whenever they will be developed in future, we can embed them within our greedy (or pseudo-greedy) algorithm to provide solutions for the PMCLP-PC-QoS along with their computational complexity.
\end{remark} \vspace{-0em}

\subsection{Approximation Ratios} \label{sec:ApproxRatio} \vspace{-0em}

In the following theorem, we provide approximation ratios associated with the greedy and pseudo-greedy algorithms for the PMCLP-PC-QoS. Let $(\bx^*, \by^*, \bz^*) \in \sr^{2p} \times \Xi$ be an optimal solution and $(\bx^a, \by^a, \bz^a) \in \sr^{2p} \times \Xi$ be an approximate solution for the PMCLP-PC-QoS. Then the approximation ratio corresponding to the approximate solution (or algorithm) is $\gamma_a = \frac{f({\bx^a, \by^a, \bz^a})}{ f({\bx^*, \by^*, \bz^*})}$.

%

\begin{theorem} \label{thm:GreedyAlgo}
Let the approximation ratios for the greedy algorithm (Algorithm 1) and the pseudo-greedy algorithm (Algorithm 2) be denoted by $\gamma_g$ and $\gamma_{r}$, respectively. Then 
\begin{equation} \label{eq:thm1}
\gamma_g > 1- \frac{1}{e} \text{ and } \gamma_{r} > 1- \frac{1}{e^{\eta}}, 
\end{equation}
where $e$ is the base of natural logarithm.
\end{theorem}

\noindent {\it Proof.} Recall that the notations $\psi^j_g$ in Algorithm 1 and $\kappa^j_{r}$ in Algorithm 2 denote the covered reward returned by the functions \texttt{SingleSZProblem{$( \mathcal{D}^g_j, s_j)$}} and  \texttt{$\eta$-ApproxSSP{$( \mathcal{D}^r_j, s_j)$}}, respectively, for $j \in \{1,\ldots,p \}$. In other words, $\psi^j_g$ is the  maximum reward and $\kappa^j_{r}$ is the $\eta$-approximate reward covered by the SZ $s_j$ for the given set of DZs $\mathcal{D}^g_j$ and $\mathcal{D}^r_j$, respectively. This implies that $\kappa^j_r \geq \eta \psi^j_r$ where $\psi^j_r$ is the maximum reward covered by the SZ $s_j$ for $\mathcal{D}^r_j$ set of DZs. For the sake of convenience, in this proof, we denote the optimal solution value for the PMCLP-PC-QoS instance, i.e., $f({\bx^*, \by^*, \bz^*})$, by $f^*$. Therefore, the approximation ratio for the greedy algorithm (Algorithm 1) and the pseudo-greedy algorithm (Algorithm 2) are given by 
\begin{equation}
\displaystyle \gamma_g = \frac{1}{f^*} \lp \sum_{l=1}^p \psi^l_g \rp \text{ and } \gamma_{r} = \frac{1}{f^*} \lp \sum_{l=1}^p \kappa^l_r \rp,
\end{equation}
respectively. Let the optimal solution value of PMCLP-PC-QoS instance with $\mathcal{D}^g_j$ (or $\mathcal{D}^r_j$) as the input set of DZs be denoted by~$\zeta^g_j$ (or $\zeta^r_j$). Then, based on Observation 1, 
\begin{align} \label{eq:obs1}
p \psi^j_g \geq \zeta^g_j \text{ and } p \kappa^j_r \geq p \eta \psi^j_r \geq \eta \zeta^r_j 
\end{align}
where $\psi^j_r$ is the covered reward returned by \texttt{SingleSZProblem}($\mathcal{D}^r_j, s_j$), for $j = 1,\ldots,p$. Note that for $j=1$, $\mathcal{D}^g_1 = \mathcal{D}^r_1 = \mathcal{D}$, and hence $\zeta^g_1 = \zeta^r_1 = f^*$. Also, since after each iteration $l \in \{1,\ldots, j-1\}$ of the greedy algorithm and pseudo-greedy algorithm, we ``trim-out'' DZs of total reward $\psi^l_g$ and $\kappa^l_r$, respectively, we get 
\begin{align} \label{eq:greedy1}
p \psi^j_g \geq \zeta^g_j \geq f^* - \sum_{l=1}^{j-1} \psi^l_g \text{ and } p \kappa_r^j \geq \eta \zeta^r_j \geq \eta \lp f^* - \sum_{l=1}^{j-1} \kappa^l_r \rp,
\end{align}
where $\psi^0_g = \kappa^0_r =0$. This implies that 
\begin{align} \label{eq:greedy2}
p \sum_{l=1}^{j} \psi^l_g & \geq f^* + \lp p-1 \rp \sum_{l=1}^{j-1} \psi^l_g \geq f^* \lp 1 + \frac{p-1}{p} + \lp\frac{p-1}{p}\rp^2 + \ldots + \lp\frac{p-1}{p} \rp^{j-1}  \rp
\end{align}
and
\begin{align} \label{eq:pseudogreedy2}
p \sum_{l=1}^{j} \kappa^l_r & \geq \eta f^* + \lp p-\eta \rp \sum_{l=1}^{j-1} \kappa^l_r \geq \eta f^* \lp 1 + \frac{p-\eta}{p} + \lp\frac{p-\eta}{p}\rp^2 + \ldots + \lp\frac{p-\eta}{p} \rp^{j-1} \rp.
\end{align}
For $j=p$, Inequalities \eqref{eq:greedy2} and \eqref{eq:pseudogreedy2} reduce to
\begin{align*}
\gamma_g = \frac{1}{f^*} \lp \sum_{l=1}^p \psi^l_g \rp \geq \lp 1 - \lp \frac{p-1}{p} \rp^p  \rp \text{ and } \gamma_r = \frac{1}{f^*} \lp \sum_{l=1}^p \kappa^l_r \rp \geq \lp 1 - \lp \frac{p-\eta}{p} \rp^p  \rp,
\end{align*}
respectively. Now because the functions in the right-hand sides of the last two inequalities are decreasing in $p$, we compute limit of these functions as $p$ approaches infinity and get 
$\gamma_g > 1- 1/e \text{ and } \gamma_{r} > 1- 1/e^{\eta}$. 

\begin{corollary}
For a given $p \geq 1$, the greedy and pseudo-greedy algorithms provide solutions whose values are at least $1- [(p-1)/p]^p$ and $1- [(p - \eta)/p]^p$, respectively, times the optimal value for the PMCLP-PC-QoS instance with $\Xi_s = \Xi$ for all $s \in \mathcal{S}$. 
\end{corollary}

\section{Exact Algorithm for PMCLP-PCR-QoS}\label{sec:PMCLP-PCR-QoS}

In this section, we analyze the objective function of PMCLP-PCR-QoS with axis-parallel rectangular DZs and SZs. We prove some theoretical properties that help us validate the reduction of search space for an optimal solution. We also present an implicit enumeration technique, a customized branch-and-bound based exact algorithm, to  solve PMCLP-PCR-QoS. Then, we computationally evaluate performance of this approach. As mentioned before, \cite{BanKia17IJOC} presented a branch-and-bound algorithm and theoretical properties of the objective function for PMCLP-PCR-QoS with $\Xi = \{1\}$, i.e., all SZs have same and fixed dimensions or QoS (denoted by PMCLP-PCR). Our algorithm not only generalizes their approach, but we further strengthen the theoretical properties provided for PMCLP-PCR, thereby improving the computational efficiency of their approach.

\subsection{Theoretical Properties for PMCLP-PCR-QoS}


We first extend some definitions introduced for PMCLP-PCR \citep{BanKia17IJOC} and introduce some new definitions.

\begin{definition}[DZ Critical Values, denoted by DCVs]\label{def:DCVs}
For each DZ $d_i$ and scaling factor $z\in \Xi$, we define a set of $x$ DCVs and a set of $y$ DCVs as follows.
\begin{align*}
X_D^{d_i,z} &:= & \big\{x_{d_i,z}^{O1}  &= x_{d_i}- w_{s_0} z, & x_{d_i,z}^{I1}  &= x_{d_i},  &x_{d_i,z}^{I2}  &=x_{d_i}+w_{d_i}- w_{s_0} z, &  x_{d_i,z}^{O2} &=x_{d_i} + w_{d_i}\big\}, \\
Y_D^{d_i,z} &:= & \big\{y_{d_i,z}^{O1}  &=y_{d_i}- l_{s_0} z,  &  y_{d_i,z}^{I1}  &=y_{d_i},  & y_{d_i,z}^{I2}  &=y_{d_i}+l_{d_i}- l_{s_0}z,   & y_{d_i,z}^{O2}  &=y_{d_i}+l_{d_i}\big\}.
\end{align*}
We denote the set of all $x$ DCVs and the set of all $y$ DCVs for a fixed scaling factor $z \in \Xi$ by $X_D^z = \bigcup_{i} X_D^{d_i,z}$ and $Y_D^z= \bigcup_{i} Y_D^{d_i,z}$, respectively. Likewise, we denote the set of all possible $x$ DCVs (or $y$ DCVs) by $X_D = \bigcup_{z\in \Xi} X_D^z$ (or $Y_D = \bigcup_{z\in \Xi} Y_D^z$).

\end{definition}

\begin{definition}[Inner and Outer DCVs, denoted by IDCVs and ODCVs, respectively]
We classify DCVs into two categories: Inner DCVs and Outer DCVs. Specifically, for a given DZ $d$ and scaling factor $z$, $x_{d,z}^{O1}$ and $x_{d,z}^{O2}$ are $x$ ODCVs and $x_{d,z}^{I1}$ and $x_{d,z}^{I2}$ are $x$ IDCVs. Let $X_{ID}^z:=\{x_{d,z}^{I1}, x_{d,z}^{I2}\}_{d\in \mathcal{D}}$ and $X_{OD}^z:=\{x_{d,z}^{O1}, x_{d,z}^{O2}\}_{d\in \mathcal{D}}$, for $z\in \Xi$. We denote the set of all possible $x$ IDCVs and $x$ ODCVs by $X_{ID} = \cup_{z\in\Xi} X_{ID}^z$ and $X_{OD} = \cup_{z\in \Xi} X_{OD}^z$, respectively. 
We also define 
\begin{align*}
 {D_{IX}^z}(x)=\bigg\{i\in N: x\in \big\{x^{I1}_{d_i,z}, x^{I2}_{d_i,z}\big\}\bigg\} \ \ \text{ and } \ \ {D_{OX}^z}(x)=\bigg\{i\in N: x\in \big\{x^{O1}_{d_i,z}, x^{O2}_{d_i,z}\big\}\bigg\}
\end{align*}
 for given $x\in \sr$ and $z \in \Xi$. Likewise, sets $Y_{ID}^z$ , $Y_{OD}^z$, $Y_{ID}$, $Y_{OD}$, ${D_{IY}^z}(y)$, and ${D_{OY}^z}(y)$ are defined for $y$ DCVs. 
\end{definition}

\begin{figure}[h]
\centering
\includegraphics[width=0.6\textwidth]{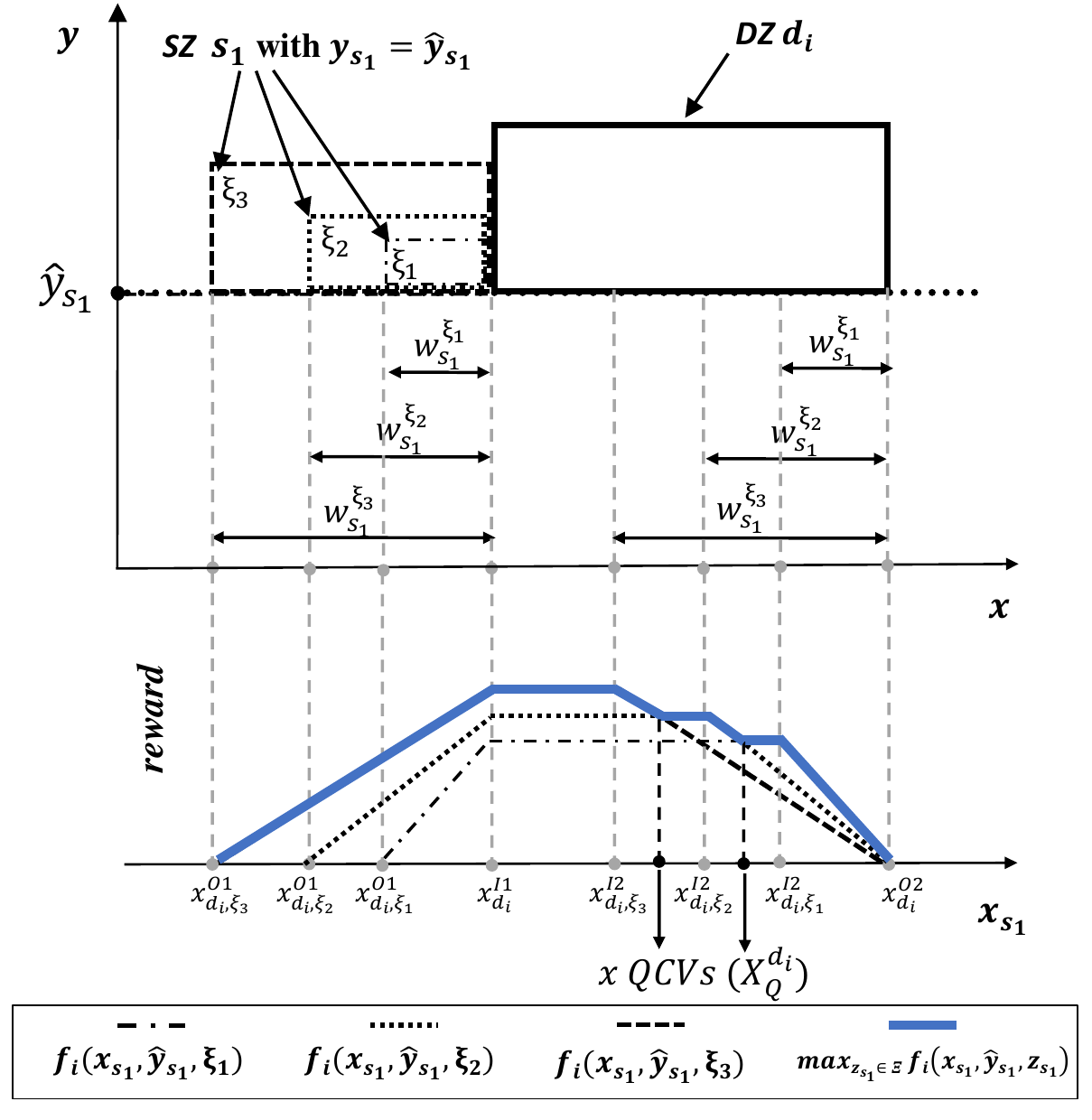} 
\caption{An example of reward function with break points at DCVs and QCVs for DZ $d_i$, SZ $s_1$ (with fixed $y_{s_1}= \hat{y}_{s_1})$ and $\Xi = \{\xi_1, \xi_2, \xi_3\}$ .} 
\label{fig:QCV}
\end{figure}

We consider a DZ $d_i$ and a SZ $s_1$ with known $y_{s_1} = \hat{y}_{s_1}$ (y coordinate of its lower left corner). Figure \ref{fig:QCV} shows that $f_i(x_{s_1},\hat{y}_{s_1},\xi_1), f_i(x_{s_1},\hat{y}_{s_1},\xi_2)$ and $f_i(x_{s_1},\hat{y}_{s_1},\xi_3)$, i.e., reward function in \eqref{eq:MaxRewardFunc} for a single SZ with fixed scaling factor $z_{s_1} \in \Xi=\{\xi_1,\xi_2,\xi_3\}$, are piecewise linear functions of $x_{s_1}$ with breakpoints at x DCVs, i.e., $X^{d_i,z_{s_1}}_D$. 
Now, note that $\max_{z_{s_1} \in \Xi} f_i(x_{s_1}, \hat{y}_{s_1}, z_{s_1})$ is also piecewise linear (represented by bold blue line segments in Figure \ref{fig:QCV}) which has breakpoints not only at x DCVs, but also at so-called $x$ QoS critical values ($x$ QCVs). A set of $x$ QCVs for DZ $d_i$, denoted by $X_{Q}^{d_i}$, includes all break points $\hat{x}_{s_1} \notin X_D$ of $\max_{z_{s_1} \in \Xi} f_i(x_{s_1}, \hat{y}_{s_1}, z_{s_1})$ where functions $f_i(x_{s_1},\hat{y}_{s_1}, \hat{z}_{s_1})$ for different QoS or $\hat{z}_{s_1}\in \Xi$ intersect. We denote the set of all $x$ QCVs by $X_{Q} = \bigcup_{_{i}}X_{Q}^{d_i}$, and define ${D_{QX}}(x)=\big\{i\in N: x\in X^{d_i}_{Q}\big\}$ for given $x \in \mathbb{R}$. Likewise, sets $Y_{Q}^{d_i}$, $Y_{Q} = \cup_{_{i}}Y_{Q}^{d_i}$ and $D_{QY}(y)$ are defined for $y$ QCVs.

\begin{observation}
Function $\max_{z_{s_1} \in \Xi}f_i(x_{s_1},\hat{y}_{s_1}, z_{s_1})$ has locally convex break point at $x_{s_1} \in X^{d_i}_Q$ because maximum of two (intersecting) affine functions is locally convex at intersecting point.  
\end{observation}

\begin{definition}[SZ Critical Values, denoted by SCVs]\label{def:SCVs}
Given a scaling factor $z \in \Xi$ and a set of positioned SZs $\{s_j\}_{j \in \mathcal{J}}$, $\mathcal{J} \subset P$, such that each SZ $s_j$, $j \in \mathcal{J}$, has known coordinates of its lower left corner ($x_{s_j},y_{s_j}$) and QoS $z_{s_j}$. 
For $j \in \mathcal{J}$ and $z \in \Xi$, we define a set of $x$ SCVs and a set of $y$ SCVs as follows. 
{\fontsize{11}{11}\selectfont
\begin{align*}
\X^{s_j,{z}}_{S} &:= & \big\{x_{s_j,z}^{O1}  &= x_{s_j}- w_{s_0}z, & x_{s_j,z}^{I1}  &= x_{s_j},  & x_{s_j,z}^{I2} &=x_{s_j}+ w_{s_0} z_{s_j} - w_{s_0} z, &  x_{s_j,z}^{O2} &=x_{s_j}+ w_{s_0} {z_{s_j}}\big\},  \\
\Y^{s_j,{z}}_{S} &:= & \big\{y_{s_j,z}^{O1}  &= y_{s_j}- l_{s_0} z, & y_{s_j,z}^{I1}  &= y_{s_j},  & y_{s_j,z}^{I2} &=y_{s_j}+ l_{s_0} {z_{s_j}} - l_{s_0} z, &  y_{s_j,z}^{O2} &=y_{s_j}+ l_{s_0} {z_{s_j}}\big\}.
\end{align*}}
We denote the set of all $x$ SCVs for a specific scaling factor $z \in \Xi$ by $\X_S^z(\mathcal{J}) = \bigcup_{j \in \mathcal{J}} \X_S^{s_j,z}$. The set of all $x$ SCVs is denoted by $\X_S(\Xi,\mathcal{J})= \bigcup_{z \in \Xi}\bigcup_{j \in \mathcal{J}} \X_S^{s_j,z}$. Likewise, the set $\Y_S^z(\mathcal{J})$ and $\Y_S(\Xi,\mathcal{J})$ are defined for $y$ SCVs. Notice that when $\Xi = \{1\}$, $x^{I1}_{s_j, z} = x^{I2}_{s_j,z}$ and $y^{I1}_{s_j, z} = y^{I2}_{s_j,z}$ because $z=z_{s_j}=1$ for all $j \in P$. 
\end{definition}

\begin{definition}[Inner and Outer SCVs, denoted by ISCVs and OSCVs, respectively]
Similar to DCVs, we also classify SCVs into two categories: Inner SCVs, i.e., $\X^{s_j,{z}}_{IS}:=\{x_{s_j,z}^{I1}, x_{s_j,z}^{I2}\}$ and $\Y^{s_j,{z}}_{IS} :=\{y_{s_j,z}^{I1}, y_{s_j,z}^{I2}\}$, and Outer SCVs, i.e., $\X^{s_j,{z}}_{OS} := \{x_{s_j,z}^{O1}, x_{s_j,z}^{O2}\}$ and $\Y^{s_j,{z}}_{OS} := \{y_{s_j,z}^{O1}, y_{s_j,z}^{O2}\}$, for $z\in \Xi$ and $j \in \mathcal{J}$.   
For a fixed scaling factor $z \in \Xi$, let  $\X_{IS}^z(\mathcal{J}) = \bigcup_{j \in \mathcal{J}} \X_{IS}^{s_j,z}$ and $\X_{OS}^z(\mathcal{J}) = \bigcup_{j \in \mathcal{J}} \X_{OS}^{s_j,z}$. Also, we define sets $\X_{IS}(\Xi,\mathcal{J})= \bigcup_{z \in \Xi}\bigcup_{j \in \mathcal{J}} \X_{IS}^{s_j,z}$ and
$\X_{OS}(\Xi,\mathcal{J})= \bigcup_{z \in \Xi}\bigcup_{j \in \mathcal{J}}$ $\X_{OS}^{s_j,z}$. Similarly, $\Y_{IS}^z$, $\Y_{OS}^z$, $\Y_{IS}$ and $\Y_{OS}$ are defined for $y$ SCVs.
\end{definition}

\begin{figure}[!h]
\centering
\includegraphics[width = 0.7\textwidth, height = 0.7\textheight]{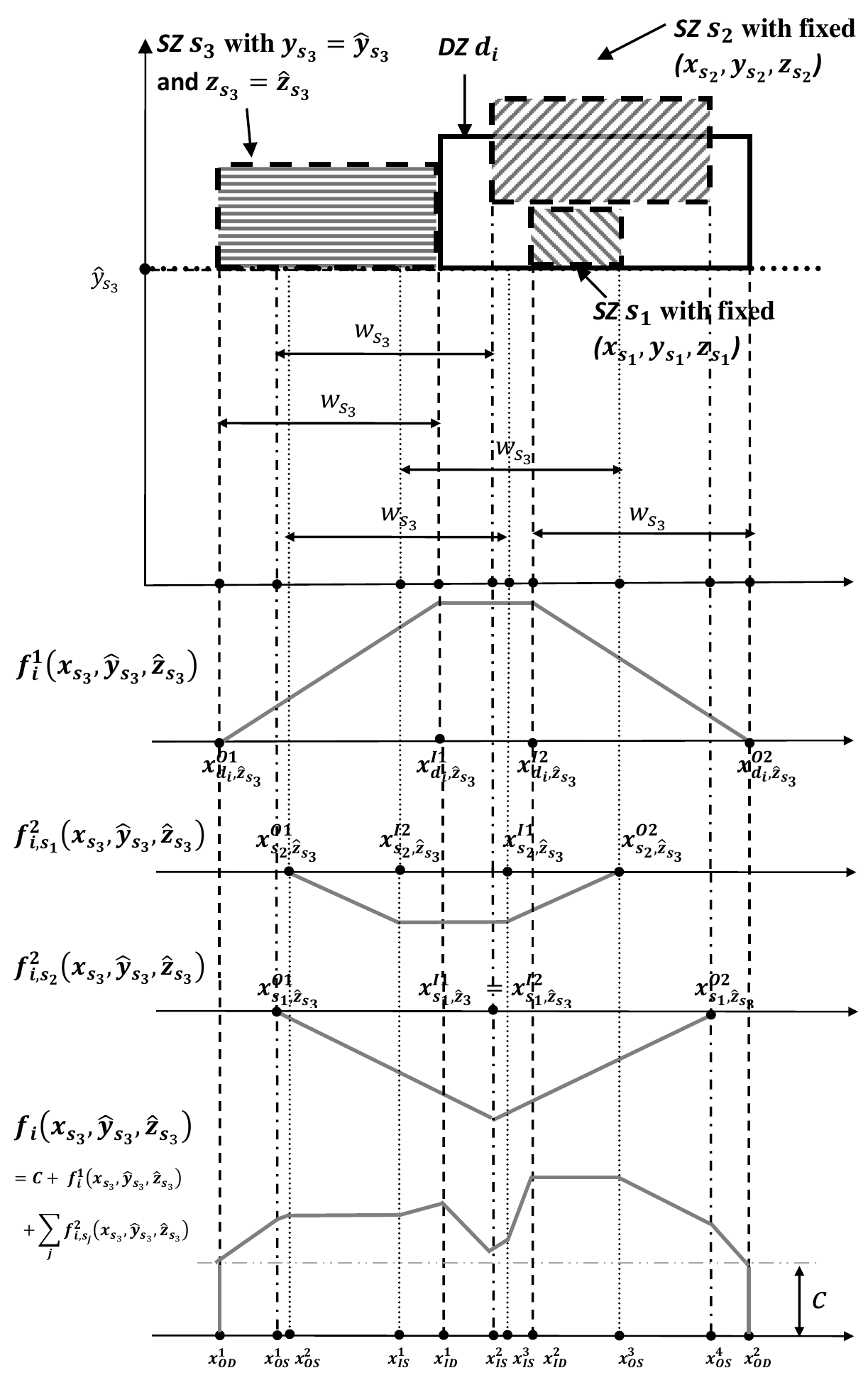} 
\caption{An example of piecewise linear reward function $f_i(\bx,\by,\bz)$ in \eqref{eq:MaxRewardFunc} with respect to argument $x_{s_3}$ (i.e., $x$ coordinate of lower left corner of SZ $s_3$) where $s_3$ has fixed scaling factor $\hat{z}_{s_3} \in \Xi$ and $y$ coordinate of its lower left corner $\hat{y}_{s_3}$, and DZ $d_i$ is partially captured by positioned SZ $s_1$ and SZ $s_2$ with fixed $(x_{s_1},y_{s_1},z_{s_1})$ and $(x_{s_2},y_{s_2},z_{s_2})$, respectively.
}
\label{fig:SCV}
\end{figure}

In Figure \ref{fig:SCV}, we provide an example to illustrate how presence of multiple SZs impact the reward function and give rise to SCVs. We consider DZ $d_i$ being partially covered by SZs $s_1$ and $s_2$ with fixed $(x_{s_1}, y_{s_1}, z_{s_1})$ and $(x_{s_2}, y_{s_2}, z_{s_2})$, respectively, such that $z_{s_1}\leq z_{s_2}$. In the rest of the paper, if a variable is not mentioned in the list of arguments that $f_i(\bx, \by, \bz)$ and $f(\bx, \by, \bz)$ take, then it means it has a fixed value. For SZ $s_3$ with $y_{s_3} = \hat{y}_{s_3}$ and $z_{s_3} = \hat{z}_{s_3} \in \Xi$, we construct the following functions: $f^1_i(x_{s_3},\hat{y}_{s_3},\hat{z}_{s_3})$ is the reward function without considering overlap of DZ $d_i$ with SZs $s_1$ and $s_2$. It is a piecewise linear function with break points at $x$ DCVs. To obtain $f_{i}(x_{s_3},\hat{y}_{s_3},\hat{z}_{s_3})$ in \eqref{eq:MaxRewardFunc} that also considers overlap of DZ $d_i$ with SZs $s_1$ and $s_2$, we define two negative reward functions $f^2_{i,{s_1}}(x_{s_3},\hat{y}_{s_3},\hat{z}_{s_3}) = - \min\{v^{\hat{z}_{s_3}}_i, v_i^{\hat{z}_{s_1}}\}  \times A(s_3 \cap s_1)$ and $f^2_{i,{s_2}}(x_{s_3},\hat{y}_{s_3},\hat{z}_{s_3})= - \min\{v^{\hat{z}_{s_3}}_i, v^{\hat{z}_{s_2}}_i\} \times A(s_3 \cap (s_2 \bs s_1))$ for SZ $s_1$ and $s_2$, respectively, then we add them to  $f^1_{i}(x_{s_3},\hat{y}_{s_3},\hat{z}_{s_3}) + c$, where $c = v_i^{\hat{z}_{s_1}} A(d_i \cap s_1) + v_i^{\hat{z}_{s_2}} A(d_i \cap (s_2 \setminus s_1))$ is the sum of reward captured by SZ $s_1$ and $s_2$ from DZ $d_i$. Observe that function $f^2_{i,{s}}(.)$, $s \in \{s_1,s_2\}$, is piecewise linear with break points at $x$ SCVs belonging to $\mathcal{X}^{s,\hat{z}_{s_3}}_S$. In case $s_1 \cap s_2 \neq \emptyset$, $f^2_{i,{s_2}}(.)$ will have break points at $x$ SCVs belonging to $\mathcal{X}^{s_1,\hat{z}_{s_3}}_S \cup \mathcal{X}^{s_2,\hat{z}_{s_3}}_S$.
As a result, function $f_{i}(x_{s_3},\hat{y}_{s_3},\hat{z}_{s_3})$ is also piecewise linear with breakpoints at $x_{s_3}$ belonging to $x$ DCVs and $x$ SCVs, i.e., $X^{d_i, \hat{z}_{s_3}}_D \cup \mathcal{X}_S^{\hat{z}_{s_3}}(\mathcal{J})$ where $\mathcal{J}=\{1,2\}$.

\begin{observation}\label{obs:ODCVnISCV}
Given a SZ $s_t \in \mathcal{S}$ and a set of positioned SZs, $\{s_j\}_{j \in \mathcal{J}}$, $\mathcal{J} = P\setminus\{t\}$ such that each SZ $s_j$, $j \in \mathcal{J}$, has known coordinates of its lower left corner and QoS, function $f_i(x_{s_t},\hat{y}_{s_t}, \hat{z}_{s_t})$ at $x_{s_t} = \hat{x}_{s_t} \in X^{\hat{z}_{s_t}}_{OD} \cup X^{\hat{z}_{s_t}}_{IS}(\mathcal{J})$ either is linear with no break point or has a locally convex break point.
\end{observation}

\begin{figure}[!h]
\centering
\includegraphics[width= 1\textwidth]{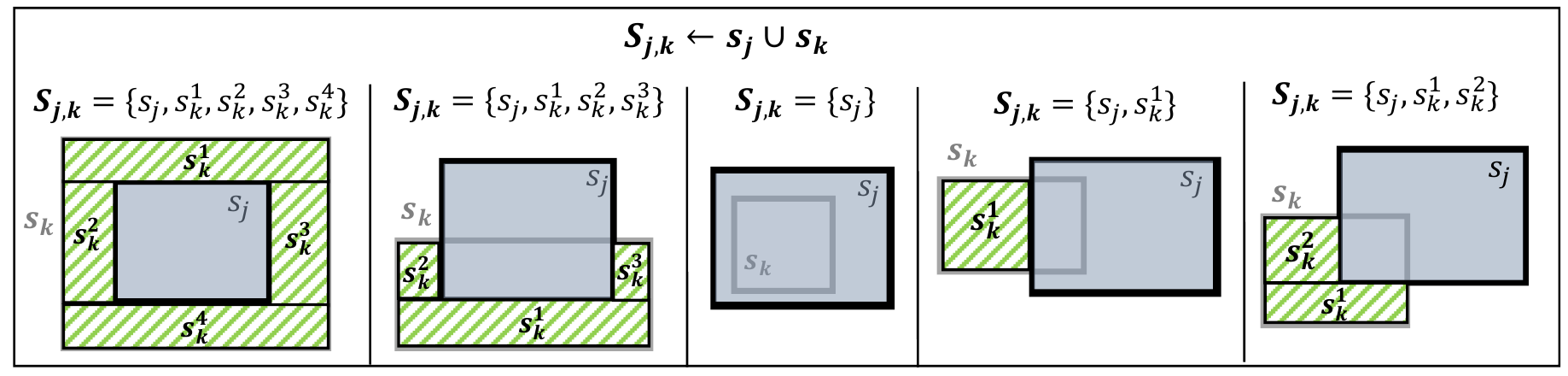}
\caption{Representing two overlapping SZs ($s_k$ and $s_j$) as a set of non overlapping rectangles}
\label{fig:fx}
\end{figure}

\begin{observation}\label{obs:gxbreakpoints}
Assume that each SZ $s \in \mathcal{S}$ has fixed coordinates of its lower left corner $(\hat{x}_s,\hat{y}_s)$ with known scaling factor $\hat{z}_s$. For each DZ $d_i$, we define a function $g_i(x) = f_i(\hat{x}_{s_1}+x , \hat{x}_{s_2}+x ,\ldots, \hat{x}_{s_p}+x, \hat{y}_{s_1},\ldots, \hat{y}_{s_p}, \hat{z}_{s_1},\ldots, \hat{z}_{s_p})$ where $x \in \sr^1$, i.e., a reward function obtained by  moving all fixed SZs simultaneously and parallel to $x$ axis. To analyze this function, we consider the following two cases. In case there is no overlap between any two SZs, i.e., $s_j \cap s_k = \emptyset$ for all $j,k \in P$ and $j\neq k$,
it is easy to observe that $g_i(x)$ is a piecewise linear function of $x$ and its break points (if exist) occur whenever $\hat{x}_{s} + x \in X^{d_i, \hat{z}_{s}}_D$ for any $s \in \mathcal{S}$. We demonstrate that these properties of $g(x)$ still hold, even in case there is overlap between SZs. Note that a set of overlapping rectangles (SZs) can be represented by non overlapping rectangles such that $x$ (or $y$) coordinate of each new vertical (or horizontal, respectively) edge, if exists, of non-overlapping rectangles coincide with $x$ (or $y$) coordinate of an edge of original overlapping SZs (see Figure \ref{fig:fx}). Consequently, the reward function obtained by moving these non overlapping rectangles simultaneously parallel to $x$ axis also has break points (if any) whenever $\hat{x}_{s} + x \in X^{d_i, \hat{z}_{s}}_D$ for any $s \in \mathcal{S}$. The same is true for a reward function obtained by simultaneously moving all SZs parallel to $y$ axis.
\end{observation}

\begin{theorem}\label{theorem:MSP-MQoS}
There exists an optimal solution $(\bx^*, \by^*, \bz^*)=(x_{s_1}^*, \ldots,x_{s_p}^*,y_{s_1}^*,\ldots,y_{s_p}^*, z_{s_1}^*, \linebreak\ldots,z_{s_p}^* )$ where $(x_{s_j}, y_{s_j}, z_{s_j}) \in \sr^2 \times \Xi$ for all $j = 1,\ldots, p$ such that the following conditions are satisfied:
\begin{description}

\item[(i)] There exists a non-repetitive sequence $s_{t_1}, s_{t_2}, \ldots, s_{t_p} \in \mathcal{S}$ of the SZs such that $x_{s_{t_1}}^* \in X_{ID}$ and for $k=2,\ldots,p$, $x_{s_{t_k}}^* \in X_{ID} \cup \mathcal{X}_{OS}(\Xi,T_k)$ where $T_k=\{t_1, \ldots, t_{k-1}\}$ and $t_k \in P$.

\item[(ii)] There exists a non-repetitive sequence $s_{q_1}, s_{q_2}, \ldots, s_{q_p} \in \mathcal{S}$ of the SZs such that $y_{s_{q_1}}^* \in Y_{ID}$ and for $k=2,\ldots,p$, $y_{s_{q_k}}^* \in Y_{ID} \cup \mathcal{Y}_{OS}(\Xi,Q_k)$ where $Q_k=\{q_1, \ldots, q_{k-1}\}$ and $q_k \in P$. %
\end{description}
In other words, $x_{s_j}^*$  is an $x$ IDCVs  or $x$ OSCVs for all $j \in P$ and $x_{s_t}^* \in X_{ID}$ for at least one SZ $s_t$, $t\in P$. The same is true for $y_{s_j}^*$, $j \in P$.
\end{theorem}

\noindent {\it Proof.} Assume there is an optimal solution $(\hat{\bx}, \hat{\by}, \hat{\bz})=(\hat{x}_{s_1}, \ldots,\hat{x}_{s_p},\hat{y}_{s_1},\ldots,\hat{y}_{s_p},\hat{z}_{s_1}, \ldots$ $,\hat{z}_{s_p})$ such that there is no $t \in P$ for which $\hat{x}_{s_t} \in X_{ID}$.
For any $t \in P$, if $\hat{x}_{s_t} \in X_{OD}$ or $\hat{x}_{s_t} \in X_Q$, then $D_{OX}(\hat{x}_{s_t}) \neq \emptyset$ or $D_{QX}(\hat{x}_{s_t}) \neq \emptyset$, respectively. Thus, for $i \in D_{OX}(\hat{x}_{s_t})$, function $f_i(x_{s_t},\hat{y}_{s_t}, \hat{z}_{s_t})$ at $x_{s_t} = \hat{x}_{s_t}$ either is linear with no break point or has a locally convex break point (Observation \ref{obs:ODCVnISCV}). Also, for $i \in D_{QX}(\hat{x}_{s_t})$, function $f_i(x_{s_t},\hat{y}_{s_t}, z_{s_t})$ has a locally convex break point at $x_{s_t} = \hat{x}_{s_t}$ and $z_{s_t}=\hat{z}_{s_t}$ (see Figure \ref{fig:QCV}).
Since there is no $t$ for which $\hat{x}_{s_t} \in X_{ID}$, then it means for all $i \in N \setminus (D_{OX}(\hat{x}_{s_t}) \cup D_{QX}(\hat{x}_{s_t}))$, $f_i(x_{s_t},\hat{y}_{s_t}, \hat{z}_{s_t})$ is linear with no break point at $x_{s_t} = \hat{x}_{s_t}$. Therefore, 
\begin{align*}
f(x_{s_t},\hat{y}_{s_t}, \hat{z}_{s_t}) = \sum_{i \in D_{OX}(\hat{x}_{s_t}) \cup D_{QX}(\hat{x}_{s_t})} f_i(x_{s_t},\hat{y}_{s_t},\hat{z}_{s_t}) + \sum_{i \in N \setminus (D_{OX}(\hat{x}_{s_t}) \cup D_{QX}(\hat{x}_{s_t}))} f_i(x_{s_t},\hat{y}_{s_t},\hat{z}_{s_t})
\end{align*}
either (a) has a locally convex break point at $x_{s_t} = \hat{x}_{s_t}$ that contradicts the optimality of $(\hat{x}_s,\hat{y}_s,\hat{z}_s )$, and hence, $\hat{x}_{s_t} \notin (X_{OD} \cup X_{Q})$, or (b) is linear with no break point which implies $\hat{x}_{s_t} \notin X_Q$. Let $g(x) = f(\hat{x}_{s_1}+x, \hat{x}_{s_2} + x, \ldots, \hat{x}_{s_p} + x, \hat{y}_{s_1},\ldots,\hat{y}_{s_p},\hat{z}_{s_1}, \ldots$ $,\hat{z}_{s_p})$ where $x\in \sr$. At $x=0$, $g(x)$ is linear and since $g(0)$ is equal to the optimal solution value, it is also constant. From Observation \ref{obs:gxbreakpoints}, we know that $g(x)$ is a piecewise linear function of $x$ and its break points occur whenever $\hat{x}_s + x \in X_D$ for any $s\in \mathcal{S}$. Let point $(\bx^*, \hat{\by}, \hat{\bz}) = (x^*_{s_1}, \ldots, x^*_{s_p},\hat{y}_{s_1},\ldots,\hat{y}_{s_p},\hat{z}_{s_1}, \ldots$ $,\hat{z}_{s_p})$ be a break point closest to $(\hat{\bx}, \hat{\by}, \hat{\bz})$, such that $x^*_{s} = \hat{x}_s + x^*$ for all $s \in \mathcal{S}$. Note that $f(\bx^*, \hat{\by}, \hat{\bz}) = f(\hat{\bx}, \hat{\by}, \hat{\bz})$. We claim that there exists $t_1\in P$ such that $x^*_{s_{t_1}} = \hat{x}_{s_{t_1}} + x^* \in X_{ID}$ because of the following reasons: \emph{Case I.} If $x^*_{s_{t_1}} \notin X_{ID}$ and $x^*_{s_{t_1}} \in X_{OD}$, then either $(\bx^*, \hat{\by}, \hat{\bz})$ is a locally convex break point which contradicts its optimality, or $g(x)$ is linear at $x^*$ with no break point. \emph{Case II.} If $x^*_{s_{t_1}} \notin X_{ID}$ and $x^*_{s_{t_1}} \in X_{Q}$, then again $(\bx^*, \hat{\by}, \hat{\bz})$ is a locally convex break point which contradicts its optimality.  The same arguments can be applied to get an optimal solution $(\bx^*, \by^*, \hat{\bz})$ such that there exists $q_1 \in P$ for which ${y}^*_{s_{q_1}} \in Y_{ID}$.

Assume that there exists an optimal solution $(\hat{\bx}, \hat{\by}, \hat{\bz})=(\hat{x}_{s_1}, \ldots,\hat{x}_{s_p},\hat{y}_{s_1},\ldots,\hat{y}_{s_p},\hat{z}_{s_1}, \ldots$ $,\hat{z}_{s_p})$ such that there exists $t_1 \in P$ for which $\hat{x}_{s_{t_1}} \in X_{ID}$, but there does not exist a non-repetitive sequence $s_{t_1}, s_{t_2}, \ldots, s_{t_p}$ of the SZs such that for $k=2,\ldots,p$, $\hat{x}_{s_{t_k}} \in X_{ID} \cup \mathcal{X}_{OS}(\Xi,T_k)$ where $t_k \in P$ and $T_k=\{t_1, \ldots, t_{k-1}\}$. This implies that for all $t \in P\bs \{t_1\}$, $\hat{x}_{s_t}$ is neither an $x$ IDCV nor an $x$ OSCV. In case $\hat{x}_{s_t} \in X^{\hat{z}_{s_t}}_{OD} \cup  \X^{\hat{z}_{s_t}}_{IS}(P\bs\{t\})$ for any $t \in P\bs\{t_1\}$, then based on Observation \ref{obs:ODCVnISCV}, $f(x_{s_t}, \hat{y}_{s_t}, \hat{z}_{s_t})$ at $x_{s_t} = \hat{x}_{s_t}$ either has a locally convex break point (which contradicts its optimality) or is linear with no break point. Likewise, in case $\hat{x}_{s_t} \in X_Q$ for any $t\in P\bs\{t_1\}$, $f(x_{s_t}, \hat{y}_{s_t}, \hat{z}_{s_t})$ at $x_{s_t} = \hat{x}_{s_t}$ has a locally convex break point, and thereby contradicts its optimality. Hence, for each $t \in P\bs\{t_1\}$, $f(x_{s_t}, \hat{y}_{s_t}, \hat{z}_{s_t})$ is linear at $x_{s_t} = \hat{x}_{s_t}$, and it is also constant at $x_{s_t} = \hat{x}_{s_t}$ because $(\hat{\bx}, \hat{\by}, \hat{\bz})$ is an optimal solution. For $t \in P\bs T_2$ where $T_2:=\{t_1\}$, let $x^*_{s_t} \in X_D \cup \X_S^{\hat{z}_{s_t}}(T_2)$ be a break point of $f(x_{s_t}, \hat{y}_{s_t}, \hat{z}_{s_t})$ that is closest to $\hat{x}_{s_t}$ and $f(x^*_{s_t}, \hat{y}_{s_t}, \hat{z}_{s_t}) = f(\hat{x}_{s_t}, \hat{y}_{s_t}, \hat{z}_{s_t})$. Now, pick $t_2 \in P\bs T_2$ such that $|\hat{x}_{s_{t_2}} - x^*_{s_{t_2}}|$ is minimum. Since $(\hat{x}_{s_{t_1}}, x^*_{s_{t_2}}, \{\hat{x}_{s_t}\}_{t \in P\bs T_3}, \hat{\by}, \hat{\bz})$ is also an optimal solution and a break point, $x^*_{s_{t_2}} \notin X_{OD} \cup \X_{IS}(\Xi, T_2) \cup X_Q$ because of the same reasons previously explained. Hence, $x^*_{s_{t_2}} \in X_{ID} \cup \X_{OS}^{\hat{z}_{s_{t_2}}}(T_2) \subseteq X_{ID} \cup \X_{OS}(\Xi, T_2)$. 
 
 Next, we use induction to show that another optimal solution $(\hat{x}_{s_{t_1}}, x^*_{s_{t_2}}, x^*_{s_{t_3}}, \ldots, x^*_{s_{t_p}}, \hat{\by}, \hat{\bz})$ can be found such that $\hat{x}_{s_{t_1}} \in X_{ID}$ and for $k=2,\ldots,p$, $x^*_{s_{t_k}} \in X_{ID} \cup \mathcal{X}_{OS}(\Xi,T_k)$ where $t_k \in P$ and $T_k=\{t_1, \ldots, t_{k-1}\}$. To do so, assume that there exists $b \in \{2,\ldots, p\}$ and an optimal solution $\lp \hat{x}_{s_{t_1}}, x^*_{s_{t_2}}, \ldots, x^*_{s_{t_b}}, \{\hat{x}_{s_t}\}_{t \in P\bs T_{b+1}}, \hat{\by}, \hat{\bz}\rp$ such that $x^*_{s_{t_k}} \in X_{ID} \cup \mathcal{X}_{OS}(\Xi,T_k)$ for $k=2,\ldots,b$, and $\hat{x}_{s_{t}} \notin X_{ID} \cup \mathcal{X}_{OS}(\Xi,T_k)$ for $t \in P \bs T_{b+1}$. Note that for $b=2$, we have already proved the existence of such optimal solution and in case $b=p$, then $P\bs T_{b+1} = \emptyset$ and hence we have an optimal solution that satisfies condition ($i$). Now for $b \in \{2, \ldots, p-1\}$, we have to find another optimal solution $\lp \hat{x}_{s_{t_1}}, x^*_{s_{t_2}}, \ldots, x^*_{s_{t_b}}, x^*_{s_{t_{b+1}}}, \{\hat{x}_{s_t}\}_{t \in P\bs T_{b+2}}, \hat{\by}, \hat{\bz}\rp$ such that $x^*_{s_{t_{b+1}}} \in X_{ID} \cup \mathcal{X}_{OS}(\Xi,T_{b+1})$ where $t_{b+1} \in P\bs T_{b+1}$. Recall that for each $t \in P\bs\{t_1\}$, $f(x_{s_t}, \hat{y}_{s_t}, \hat{z}_{s_t})$ is linear at $x_{s_t} = \hat{x}_{s_t}$, and it is also constant at $x_{s_t} = \hat{x}_{s_t}$. Since $P\bs T_{b+1} \subseteq P\bs \{t_1\}$, for $t \in P\bs T_{b+1}$, $f(x_{s_t}, \hat{y}_{s_t}, \hat{z}_{s_t})$ is linear with no break point at $\lp \hat{x}_{s_{t_1}}, x^*_{s_{t_2}}, \ldots, x^*_{s_{t_b}}, \{\hat{x}_{s_t}\}_{t \in P\bs T_{b+1}}, \hat{\by}, \hat{\bz}\rp$. Let $x^*_{s_t} \in X_D \cup \X_S^{\hat{z}_{s_t}}(T_{b+1})$, $t \in P\bs T_{b+1}$, be a break point of $f(x_{s_t}, \hat{y}_{s_t}, \hat{z}_{s_t})$ that is closest to $\hat{x}_{s_t}$
 and $f(x^*_{s_t}, \hat{y}_{s_t}, \hat{z}_{s_t}) = f(\hat{x}_{s_t}, \hat{y}_{s_t}, \hat{z}_{s_t})$. Select index $t_{b+1} \in P\bs T_{b+1}$ such that $\big|\hat{x}_{s_{t_{b+1}}} -   x^*_{s_{t_{b+1}}}\big|$ is minimum to get another optimal solution $\lp \hat{x}_{s_{t_1}}, x^*_{s_{t_2}},\ldots,x^*_{s_{t_b}}, x^*_{s_{t_{b+1}}} \{\hat{x}_{s_t}\}_{t \in P\bs T_{b+2}}, \hat{\by}, \hat{\bz}\rp$ where $x^*_{s_{t_{b+1}}} \in X_{ID} \cup \X_{OS}^{\hat{z}_{s_{t_{b+1}}}}(T_{b+1}) \subseteq X_{ID} \cup \X_{OS}(\Xi, T_{b+1})$. Note that $x^*_{s_{t_{b+1}}} \notin X_{OD} \cup \X_{IS}(\Xi, T_{b+1}) \cup X_Q$ because otherwise $f(x_{s_{t_{b+1}}}, \hat{y}_{s_{t_{b+1}}}, \hat{z}_{s_{t_{b+1}}})$ either has a locally convex break point or is linear with no break point at $x_{s_{t_{b+1}}} = x^*_{s_{t_{b+1}}}$. 
 This implies that by sequentially setting $b=2, \ldots, p-1$, we can eventually get an optimal solution that satisfies condition $(i)$. 
 %
 By applying the same arguments for $y$ coordinates,  we reach an optimal solution that satisfies condition~($ii$). \qed

 \begin{remark}
For $m=1$, Theorem \ref{theorem:MSP-MQoS} in this paper strengthens Theorem 3 of \cite{BanKia17IJOC}. 
According to the latter, there exists an optimal solution for PMCLP-PCR-QoS with $\Xi = \{1\}$ (PMCLP-PCR), denoted by $(x^*_{s_1}, \ldots, x^*_{s_p}, y^*_{s_1}, \ldots, y^*_{s_p})$, such that $x^*_{s}$ is an $x$ IDCV or $x$ SCV for all $j \in P$. Whereas, according to Theorem \ref{theorem:MSP-MQoS}, there exists an optimal solution for PMCLP-PCR such that $x^*_{s}$ is an $x$ IDCV or $x$ OSCV for all $j \in P$, thereby reducing the solution search space. In Section \ref{sec:ComExp}, we will observe how this improves the exact algorithm of~\cite{BanKia17IJOC} for PMCLP-PCR and reduces the time taken to solve it.
\end{remark}

\subsection{Exact Branch-and-Bound Algorithm for PMCLP-PCR-QoS}

In this section, we present a generalization of the branch-and-bound (B\&B) algorithm presented by \cite{BanKia17IJOC} for PMCLP-PCR-QoS with $\Xi=\{1\}$ to solve PMCLP-PCR-QoS with $\Xi = \{\xi_1,\xi_2,\ldots, \xi_m\}$ where $m \geq 1$, thereby allowing adjustable QoS for each SZ. The former can be utilized to solve PMCLP-PCR-QoS with an additional restriction of $z_{s_j} = z \in \Xi$ (same QoS for all SZs) by solving PMCLP-PCR $m$ number of times. For $m=1$, our new algorithm implicitly enumerates a reduced solution search space, in comparison to the solution space considered in \cite{BanKia17IJOC}, and as a result, our proposed algorithm is computationally faster (see Section \ref{sec:ComExp} for more details). For $m\geq 2$, though the outline of our approach is similar to a generic B\&B approach, its key components such as branching, upper bound computation, lower bound updates, and node selection, as well as the data structures to efficiently implement them are customized for PMCLP-PCR-QoS. It is important to note that PMCLP-PCR-QoS can also be solved by brute-force search, i.e., explicit enumeration of all possible combinations of $x$ and $y$ critical values, in particular IDCVs and OSCVs. However, this approach is computationally very expensive. For an example, to solve an instance of PMCLP-PCR-QoS with $p=m=2$ and $n=100$ using explicit enumeration, we have to evaluate reward function for at least $2.56$ million solutions. In comparison, our approach explores only 13,092 nodes (average for 10 instances) of B\&B tree where for most of these nodes, it computes upper bound that is computationally less expensive than evaluating reward function. 

Since some attributes of our B\&B approach for PMCLP-PCR-QoS are inherited from B\&B approach for PMCLP-PCR, we use similar notations and nomenclatures, to the extent possible, in the ensuing subsections for the sake of reader's convenience. Before we present the algorithm, we highlight its main additional features. In PMCLP-PCR-QoS, there is an additional set of decisions of choosing QoS (i.e., scaling factor) for each SZ. Therefore, we need different strategies and data structures to incorporate these decisions in building B\&B tree, designing branching routine, and computing upper bound at each node (see Sections \ref{sec:Branching} and \ref{sec:UpperBound}, respectively, for more details). The branching strategies also influence the generation of OSCVs that depends on scaling factor (dimensions) of SZs. Moreover, to compute lower bound at a leaf node (covered reward by $p$ SZs having different QoS), we ensure that for a subregion (partially) covered by multiple SZs with different QoS, the reward rate corresponding to the best QoS among these SZs is considered in the reward calculations. 
In B\&B approach for PMCLP-PCR, the authors construct two auxiliary trees (one for each of $x$ and $y$ axes) to store partitions of initial set of IDCVs. This construction provides computational advantages and avoid repetition of computing partitions of the set of IDCVs for different SZs. However, in order to extend this procedure for PMCLP-PCR-QoS, we would need $2m$ number of such auxiliary trees (two for each scaling factor). Therefore, instead of using the auxiliary trees, we evaluate four different procedures to design the B\&B tree, and consider the one that helps us to focus on the regions with high reward, provides a strong lower bound quickly, and makes the overall algorithm computationally efficient (refer to Section \ref{sec:Branching} for more details).

\subsubsection{Main Body}

Algorithm \ref{Algo:MainBody} provides a pseudocode for the proposed customized B\&B algorithm. We initialize the algorithm (Line \ref{Alg:MB:init}) by computing an initial feasible solution and lower bound of the objective function of PMCLP-PCR-QoS using \fun{GreedyAlgorithm} for the set of DZs $\mathcal{D}$ and $p$ SZs (discussed in Section \ref{sec:Greedy}). Recall that according to Theorem \ref{thm:GreedyAlgo}, this lower bound is at least $63.2\%$ of the optimal solution value for PMCLP-PCR-QoS instance. Thereafter, we construct the root node $Q_0$ of the B\&B tree in line \ref{Alg:MB:Q0} that stores $p$ sets of all $x$ IDCVs, $X^{s_j}_{Q_0} = X_{ID} = \cup_{_{z \in \Xi}} X^{z}_{ID}$, $j \in P$, and $y$ IDCVs, $Y^{s_j}_{Q_0} = Y_{ID} = \cup_{_{z \in \Xi}} Y^{z}_{ID}$ for all $j \in P$. Similarly, at each node $Q$ of the B\&B tree, we store $p$ subsets of $x$ CVs (IDCVs and OSCVs), denoted by $X^{s_j}_Q, j \in P$, and $p$ subsets of $y$ CVs (IDCVs and OSCVs), denoted by $Y^{s_j}_Q,j\in P$, along with three indicators: $(i)$ a $p\times 1$ vector $\mathcal{Z}(Q)$ to store scaling factor associated with each SZ, $(ii)$ $\mathcal{BA}(Q)$ to store branching axis, i.e $x$ or $y$, and $(iii)$~$\mathcal{BS}(Q) \in \mathcal{S}$ to store branching SZ and $index(\mathcal{BS}(Q))\in P$ denotes the index of the SZ. At root node $Q_0$, we set $\mathcal{Z}(Q_0) \leftarrow {\bf 0}$ because scaling factor is not yet selected for each SZ, $\mathcal{BA}(Q_0) \leftarrow ``X"$, $\mathcal{BS}(Q_0) \leftarrow s_1$, and $index(\mathcal{BS}(Q_0)) \leftarrow 1$. Let \texttt{LCN} denotes a list of candidate nodes that is initialized by including $Q_0$ in it. Each node in this list represents a restricted PMCLP-PCR-QoS problem (or subproblem) where SZ $s_j$, $j \in P$, has $(x_{s_j}, y_{s_j}) \in X^{s_j}_Q \times Y^{s_j}_Q$ and $z_{s_j}$ is the $j^{\text{th}}$ element of $\mathcal{Z}(Q)$. A zero element of $\mathcal{Z}(Q)$ implies that the scaling factor of the corresponding SZ has not been selected yet and it belongs to $\Xi$.

\begin{algorithm}[!h]
  \fontsize{9}{12}\selectfont
  \caption{Exact Branch-and-Bound Algorithm: Main Body} \label{Algo:MainBody}
\begin{algorithmic}[1]
\State Using \fun{GreedyAlgorithm($\cD, \Xi, p$)}, obtain an initial feasible solution and lower bound (LB);
\label{Alg:MB:init}
\State Construct a root node $Q_0$ using $x$ and $y$ IDCVs with $\mathcal{Z}(Q_0) \leftarrow {\bf 0}$, $\mathcal{BA}(Q_0) \leftarrow X$, and $\mathcal{BS}(Q_0) \leftarrow 1$;
\label{Alg:MB:Q0}
\State Initialize a list of candidate nodes $\texttt{LCN}\leftarrow \{Q_0\}$;

\State Select a node $Q$ belonging to \texttt{LCN} using depth first strategy;
\label{Alg:MB:LCN}
\State Calculate an upper bound for node $Q$:
\label{Alg:MB:UBcalc}
\Statex \hspace{2.0cm} $UB(Q) \leftarrow \fun{UpperBound($Q, \cD,\Xi,p$)}$;
\If{$UB(Q) \leq LB$} {prune the node $Q$ and remove it from \texttt{LCN}}
\label{Alg:MB:UB<LB prune}
\ElsIf{Node $Q$ is a leaf node}{\hspace{0ex} $LB \leftarrow UB(Q)$ and remove $Q$ from \texttt{LCN}}
\label{Alg:MB:leafnode prune}
\label{Alg:MB:feasiblesol}
\label{Alg:MB:LBupdate}
\Else{ Create child nodes of node $Q$ using \fun{Branching($Q$)} routine and add them to \texttt{LCN};}
\label{Alg:MB:callBranch}
\EndIf
\State Go to Step \ref{Alg:MB:LCN} if \texttt{LCN} is nonempty;\label{Alg:MB:WhileEnd} 
\State \Return{$LB$ and Optimal solution $(x_{s_1}^*, \ldots, x_{s_p}^*, y_{s_1}^*, \ldots, y_{s_p}^*,z_{s_1}^*,\ldots,z_{s_p}^* )$};
\label{Alg:MB:return}
	\end{algorithmic}
\end{algorithm}

To implicitly enumerate the solution space, we select a node $Q$ belonging to \texttt{LCN} and calculate an upper bound ($\UB$) of the optimal objective value for the restricted problem associated to node $Q$, using \fun{UpperBound} function. Refer to Section \ref{sec:UpperBound} for pseudocode and more details regarding this function. To avoid explicit enumeration of the critical points, a strong upper bound is needed so that nodes (or subregions) that cannot provide a feasible solution better than the best known solution, be pruned. In line \ref{Alg:MB:UB<LB prune}, we prune a node if $\UB \leq \LB$ (best known lower bound) and remove it from \texttt{LCN}. Otherwise, if node $Q$ is a leaf node, i.e., $|X^{s_j}_Q|=|Y^{s_j}_Q|=1$ for all $j \in P$, and $\UB > \LB$, the \fun{UpperBound} function returns an improved lower bound for the original problem (Line \ref{Alg:MB:LBupdate}), and the leaf node is pruned and removed from \texttt{LCN}. However, in case node $Q$ is not a leaf node and $\UB > \LB$, we create child nodes of node $Q$ (also referred to as parent node) using \fun{Branching} function and add them to \texttt{LCN} (Line \ref{Alg:MB:callBranch}). A child node is a subproblem that inherits (some) restrictions from its parent node along with new restrictions; refer to Section \ref{sec:Branching} for pseudocode and more details of the \fun{Branching} function. We use depth-first strategy to traverse the B\&B tree. This implies that in case the current node is pruned, we select a node that is the right sibling of the node or one of its ancestor, and in case the current node is branched, we select its first (or left most child) from \texttt{LCN} and repeat Lines \ref{Alg:MB:LCN}-\ref{Alg:MB:WhileEnd}. The algorithm is terminated when \texttt{LCN} is empty and it returns an optimal solution, i.e., location and scaling factor of each SZ $(x_{s_1}^*, \ldots, x_{s_p}^*, y_{s_1}^*, \ldots, y_{s_p}^*,z_{s_1}^*,\ldots,z_{s_p}^* )$ as well as the optimal captured reward value (i.e., best known $LB$).

\subsubsection{Branching} \label{sec:Branching}

This function creates a list of child nodes for a given parent node $Q$ (refer to Algorithm~\ref{Algo:Branch} for psuedocode). These child nodes are added in \texttt{LCN} and the most left child node of $Q$ is selected from \texttt{LCN} in the next iteration. At any given node $Q$, indicators $\mathcal{BA}(Q)$ and $\mathcal{BS}(Q)$ determine branching axis and branching SZ. These two properties are defined to select a set among the sets $X^{s_j}_Q$ and $Y^{s_j}_Q$ for $j \in P$ that is to be partitioned to create new child nodes. We use the function \funi{NewChildNode(Q)} to create child node $T$ for node $Q$, which inherits the properties of its parent. We classify the branching procedure into three categories:

\noindent (a) \emph{QoS Branching}. When $X^{\mathcal{BS}(Q)}_Q = X_{ID}$ (or $Y^{\mathcal{BS}(Q)}_Q = Y_{ID}$), i.e., the scaling factor for SZ $\mathcal{BS}(Q)$ has not been selected, we create $|\Xi|=m$ child nodes for a parent node $Q$. Each child node corresponds to a fixed scaling factor $z \in \Xi=\{\xi_1, \ldots, \xi_m\}$ for the branching SZ. As a result, a child node $T$ of node $Q$ corresponding to scaling factor $z \in \Xi$ has $\mathcal{Z}(T, \mathcal{BS}(T)) = z$ and subsets $X^{\mathcal{BS}(T)}_T= X^{Z(T, {\mathcal{BS}(T)})}_{ID}$ and $Y^{\mathcal{BS}(T)}_T= Y^{Z(T, \mathcal{BS}(T))}_{ID}$ (Lines \ref{Alg:QoSBranchStart}-\ref{Alg:QoSBranchEnd} of Algorithm~\ref{Algo:Branch}). Note that $\mathcal{Z}(T,s_j)$ denotes the $j^{th}$ element of vector $\mathcal{Z}(T)$.

\noindent (b) \emph{Partition-Based Branching}. In case $\mathcal{BA}(Q)=``X"$ (or $``Y"$), $X^{\mathcal{BS}(Q)}_Q \neq X_{ID}$ (or $Y^{\mathcal{BS}(Q)}_Q \neq Y_{ID}$), and $|X^{\mathcal{BS}(Q)}_Q| \neq 1$ (or $|Y^{\mathcal{BS}(Q)}_Q| \neq 1$), node $Q$ is branched to child nodes by partitioning the set $X^{\mathcal{BS}(Q)}_Q$ (or $Y^{\mathcal{BS}(Q)}_Q$). Each child node $T$ of parent node $Q$ corresponds to a partition where the partition is assigned to $X^{\mathcal{BS}(T)}_T$ (or $Y^{\mathcal{BS}(T)}_T$) of node $T$ (lines \ref{Alg: xpartStart}-\ref{Alg: xpartEnd}). In addition, whenever $X^{\mathcal{BS}(Q)}_Q= X^{Z(Q,\mathcal{BS}(Q))}_{ID}$ (or $Y^{\mathcal{BS}(Q)}_Q= Y^{Z(Q,\mathcal{BS}(Q))}_{ID}$) and at least one SZ, other than $\mathcal{BS}(Q)$, has fixed $x$ (or $y$) coordinates, we also create child node $T$ corresponding to each $x$ (or $y$) OSCVs generated by SZs with fixed $x$ (or $y$) positions, where SZ $\mathcal{BS}(Q)$ can be placed (lines \ref{Alg: xOSCVStart}-\ref{Alg: xOSCVEnd}). Each child node $T$ defines a further restricted problem (or subproblem) in which SZ $s_j$, $j \in P$, can only be positioned at the locations $(x,y)\in X^{s_j}_T \times Y^{s_j}_T$. At node $T$, we set the indicators $\mathcal{BA}$, $\mathcal{BS}$, and $\mathcal{Z}$ same as the parent node $Q$, except when $X^{\mathcal{BS}(Q)}_T$ (or $Y^{\mathcal{BS}(Q)}_T$) becomes singleton set, i.e., SZ $\mathcal{BS}(Q)$ has a fixed $x$ (or $y$) coordinate. For the foregoing scenario, we update $\mathcal{BS}(T)$ and/or $\mathcal{BA}(T)$ so that while branching on node $T$, if needed, an appropriate set among the sets $X^{s_j}_T$ and $Y^{s_j}_T$ for $j \in P$ is selected for partitioning. It is important to note that the efficiency and effectiveness of the B\&B algorithm to solve PMCLP-PCR-QoS depends on effective generation of partitions and updates of indicators $\mathcal{BS}(T)$ and/or $\mathcal{BA}(T)$.

\begin{algorithm}[!h]
\fontsize{9}{10}\selectfont
\caption{Branching}\label{Algo:Branch}
\begin{algorithmic}[1]
\Function{\fun{Branching}}{\pmb{$Q$}}
\If{$\mathcal{BA}(Q)=$ \funi{``X"}} 
		\If{ $index(\mathcal{BS}(Q)) \leq p$ (node $Q$ has at least a SZ whose x coordinate is not fixed)}
		        \If{$X^{\mathcal{BS}(Q)}_Q$ = $X_{ID}$}
		        \label{Alg:QoSBranchStart}
		               \For{$z \leftarrow \xi_1$ to $\xi_m$}
		                    \State $T\leftarrow$ {$NewChildNode(Q)$}
		                    \State $X^{\mathcal{BS}(T)}_T\leftarrow$ {$X_{ID}^z$};  $Y^{\mathcal{BS}(T)}_T\leftarrow$ {$Y_{ID}^z$}
		                    \State $Z(T,{\mathcal{BS}(T)})\leftarrow$ {$z$}
		                \EndFor
		               \label{Alg:QoSBranchEnd}
		            \Else
		                \For{$n_x \leftarrow 1$ to $NumPartitions(X^{\mathcal{BS}(Q)}_Q)$}
		              \label{Alg: xpartStart}
		                    \State $T\leftarrow$ {$NewChildNode(Q)$}
		                    \State $X_T^{\mathcal{BS}(T)} \leftarrow$ {$GetPartition(X_Q^{\mathcal{BS}(Q)},n_x)$}
		      \label{Alg: xgetpart} 
		                    \If{$|X_T^{\mathcal{BS}(T)}|=1$}
    		                    \State$index(\mathcal{BS}(T)) \leftarrow$ {$index(\mathcal{BS}(Q)) +1$}
    		  \label{Alg: BSupdate1}               
		                    \EndIf
		                  \EndFor
		              \label{Alg: xpartEnd}
		                  
		                  \If{$X^{\mathcal{BS}(Q)}_Q = X_{ID}^{Z(Q,{\mathcal{BS}(Q)})}$}
		              \label{Alg: xOSCVStart}
		                  \State $L_x \leftarrow \{j \in P: |X^{s_j}_Q| =1\}$
		                  \For{$\hat{x} \in (\bigcup_{j\in L_x} {\X^{s_j,{Z(Q,{\mathcal{BS}(Q)})}}_{OS}}) \setminus X_{ID} $}
    		                  \State $T\leftarrow$ {$NewChildNode(Q)$}
        	                  \State $X_T^{\mathcal{BS}(T)} \leftarrow \hat{x}$   
            	               \State $index(\mathcal{BS}(T)) \leftarrow index(\mathcal{BS}(Q)) +1$
	       \label{Alg: BSupdate2} 
	       
	       \EndFor 
                        \EndIf
                        \label{Alg: xOSCVEnd}
		          \EndIf
		          \Else
		      \State $ l_1 \leftarrow$ smallest SZ index such that $x$ position of the SZ is fixed at an x IDCV \Statex \hspace{2.6cm}and y position of the SZ is not fixed;
		      \label{Alg: PermutStart}
    	      \State $\bar{L} \leftarrow $ set of SZ indices whose x position is fixed at $x$ OSCV and y position is not fixed;
    	      \label{Alg: PermutStart2}
    	      \For{$l \in \bar{L} \cup \{l_1\}$}
    	      \label{Alg: PermutStart-detail1}
    							\State $T\leftarrow$ {$NewChildNode(Q)$} 					
    	      			\State $index(\mathcal{BS}(T)) \leftarrow$ {$l$}
    	      			\State $\mathcal{BA}(T) \leftarrow$ {"Y"}
    	      \EndFor
    	      \label{Alg: PermutStart-detail2}
    	      
    		\EndIf
    		\label{Alg: PermutEnd}
    		\ElsIf{ $\mathcal{BA}(Q)=$ \funi{"Y"} } 
    	   \If{$|Y_Q^{\mathcal{BS}(Q)}|=1$}{\hspace{1ex}do Steps \ref{Alg: PermutStart} to \ref{Alg: PermutEnd}}
    	   \label{Alg: FixYBranch}
    	       \Else
    	       \For{$n_y \leftarrow 1$ to $NumPartitions(Y^{\mathcal{BS}(Q)}_Q)$}
		       \label{Alg: ypartStart} 
		                    \State $T\leftarrow$ {$NewChildNode(Q)$}
		                    \State $Y_T^{\mathcal{BS}(T)} \leftarrow$ {$GetPartition(Y_Q^{\mathcal{BS}(Q)},n_y)$}
		      \label{Alg: ygetpart}           
		                    
		      \EndFor
    	       \If{$Y^{\mathcal{BS}(Q)}_Q = Y_{ID}^{Z(Q,{\mathcal{BS}(Q)}})$}
		                  \State $L_y \leftarrow \{j \in P: |Y^{s_j}_Q| =1\}$
		                  \For{$\hat{y} \in (\bigcup_{j\in L_y} {\Y^{s_j,{Z(Q,{\mathcal{BS}(Q)})}}_{OS}}) \setminus Y_{ID} $}
    		                  \State $T\leftarrow$ {$NewChildNode(Q)$}
        	                  \State $Y_T^{\mathcal{BS}(T)} \leftarrow \hat{y}$   
	                    \EndFor 
                        \EndIf
    	       \EndIf
                \EndIf	
\EndFunction
	\end{algorithmic}
\end{algorithm}

For partitioning, we assign a priority score to each element of the set $X^z_{ID}$, $z \in \Xi$, i.e., $x_k \in X^z_{ID}$ for $k\in \{1,\ldots,|X^z_{ID}|\}$ has priority score $p(x_k) = \sum_{d \in \mathcal{D}_k} v^z_d$ where $\mathcal{D}_k:= \{d \in \mathcal{D}: x_{d} \leq x_k<x_{d}+w_{d}\}$, i.e., sum of reward rate corresponding to scaling factor $z$ for DZs such that $x_k$ lies between $x$ coordinate of vertical edges of DZ $d$. We also considered $p(x_k) = \frac{\sum_{d \in \mathcal{D}_k} v^z_d}{|\mathcal{D}_k|}$ (average of reward rates), but observed that this choice of priority score is not effective. Then, we sort all elements of the set $X^z_{ID}$ based on their priority score in the descending order. The sorted set is also denoted by $\{x_1, x_2, \ldots, x_{|X^z_{ID}|}\}$ for the sake of convenience. Similarly, the priority score $p(y_k)$ is defined for each $y_k \in Y^z_{ID}$ and elements of set $Y^z_{ID}$ are also sorted in descending order of their priority score. 
At node $Q$, we split $X^{{\mathcal{BS}(Q)}}_Q$ (or $Y^{{\mathcal{BS}(Q)}}_Q$)  by creating disjoint subsets when  $x_{i+1}-x_i > \beta \max_{j=1,\ldots,|X^{{\mathcal{BS}(Q)}}_Q|-1}\{x_{j+1} - x_j\}$ for any $x_i, i\in\{1,\ldots,|X^{{\mathcal{BS}(Q)}}_Q|-1\}$ where $\beta \in (0,1)$. These partitions are assigned to child nodes using  $GetPartition$ in Lines \ref{Alg: xgetpart} and \ref{Alg: ygetpart}. The smaller values of $\beta$ increase the number of partitions ($NumPartitions$ in Lines \ref{Alg: xpartStart} and \ref{Alg: ypartStart}), and the bigger values of $\beta$ lead to either no partitioning or very small number of  subsets, thereby reducing the computational efficiency of the B\&B algorithm. Based on our computational experiments, we select $\beta = 0.5$ as the most appropriate value resulting in an efficient algorithm. This ensures that the B\&B algorithm traverses regions with high concentration of rewards faster, and hence, it quickly finds good lower bounds. 

In our computational studies, we also considered the strategy of sorting child nodes based on their priority score in each branching step (similar to \cite{BanKia17IJOC}). Specifically, instead of first sorting the set of IDCVs based on priority scores and then creating partitions (as discussed above), the set $X^{{\mathcal{BS}(Q)}}_Q$ (or $Y^{{\mathcal{BS}(Q)}}_Q$) sorted in ascending order of IDCVs is partitioned. Then, a priority score is assigned to each child node $T$ of node $Q$, defined by ${\sum_{x_k \in \mathcal X_{T}^{\mathcal{BS}(T)}} p(x_k)}/{|X_{T}^{\mathcal{BS}(T)}|}$, and the child nodes are sorted based on their priority score. A child node with highest priority score is placed as the left-most child node of $Q$ in the B\&B tree. Though this approach is computationally comparable to our approach in terms of solution time of the B\&B algorithm, it resulted in more number of traversed nodes. We also defined priority score for child node $T$ as $\sum_{x_k \in \mathcal X_{T}^{\mathcal{BS}(T)}} p(x_k)$, but this increased the solution time. 

We update indicator $\mathcal{BS}(T)$ by incrementing $index(\mathcal{BS}(T))$ by one (Line \ref{Alg: BSupdate1} and \ref{Alg: BSupdate2}), when $\mathcal{BA}(Q) = ``X"$, $X^{\mathcal{BS}(Q)}_T$ becomes a singleton set, and $Y^j_Q = Y_{ID}$ for all $j \in P$. As a result, the sequence $s_{t_1}, s_{t_2}, \ldots, s_{t_p}$ in Theorem \ref{theorem:MSP-MQoS} is equal to $1,\ldots, p$. This is justified because each SZ has the same set of scaling factors $\Xi =\{\xi_1, \ldots, \xi_m\}$ as well as the same set of initial $x$ IDCVs  associated with it (i.e., all SZs are homogeneous at the beginning). Therefore, different sequences $s_{t_1}, s_{t_2}, \ldots, s_{t_p}$ will lead to the same set of solutions. In contrast, when all SZs have fixed $x$ coordinate, they are no longer homogeneous. So, we cannot consider $q_1, q_2, \ldots, q_p$ in Theorem \ref{theorem:MSP-MQoS} to be equal to $1,\ldots, p$, and hence a permutation branching is required. Note that we can also set $\mathcal{BA}(Q_0) = ``Y"$ at the root node $Q_0$. In that case, the sequence $q_1, q_2, \ldots, q_p$ will be equal to $1,\ldots,p$ and permutation branching will be required for branching along $x$ axis.

\noindent(c) \emph{Permutation Branching} (Lines \ref{Alg: PermutStart}-\ref{Alg: PermutEnd}). For a node $Q$ where all $p$ SZs have fixed $x$ coordinates, i.e., $\mathcal{BA}(Q) = ``X"$ and {$index(\mathcal{BS}(Q)) > p$}, and fixed QoS, i.e., $\mathcal{Z}(Q,s_j)\neq 0$ for all $j \in P$, we start the branching along the $y$ axis using $Y^j_Q$, $j \in P$. At this node, SZs are not necessarily homogeneous because they can have different $x$ coordinates and QoS. Therefore, according to Theorem \ref{theorem:MSP-MQoS}, we consider different sequences $q_1, q_2, \ldots, q_p$ of branching SZs along $y$ axis. However, we observe that because of some symmetries in the B\&B tree, we can ignore a few sequences. For an example, when all SZs are fixed on $x$ IDCVs (say $\hat{x}_{s_1}, \hat{x}_{s_2}, \ldots, \hat{x}_{s_p}$) with QoS (say $\hat{z}_{s_1}, \hat{z}_{s_2}, \ldots, \hat{z}_{s_p}$), we only consider sequence $1,\ldots, p$ for unique solutions. This is because for each sequence $q_1, q_2, \ldots, q_p$, there exists another node in the B\&B tree where all SZs are fixed on $x$ IDCVs ($\hat{x}_{s_{q_1}}, \hat{x}_{s_{q_2}}, \ldots, \hat{x}_{s_{q_p}}$) with QoS $\hat{z}_{s_{q_1}}, \hat{z}_{s_{q_2}}, \ldots, \hat{z}_{s_{q_p}}$ as the initial set of IDCVs $X_{ID}$ and scaling factors $\Xi$ are same for all SZs. In contrast, when a SZ has an $x$ coordinate fixed at an $x$ OSCV, the foregoing property does not hold. Let $\bar{L}$ be the set of SZ indices whose $x$ coordinate is fixed at $x$ OSCVs (line \ref{Alg: PermutStart2}). Then, we consider all permutations of $s_1, s_2, \ldots, s_p$ where SZs belonging to the set $\{s_j: j \in P\setminus \bar{L}\}$ are considered as identical objects. This is done in Algorithm \ref{Algo:Branch} as follows. Let $l_1$ be the smallest SZ index having $x$ positioned at $x$ IDCVs and $y$ is not fixed (line \ref{Alg: PermutStart}). Now, for each $l \in l_1 \cup \bar{L}$, we create a child node $T$ of node $Q$ where {$index(\mathcal{BS}(T)) = l$} and $\mathcal{BA}(T)=``Y"$ (line \ref{Alg: PermutStart-detail1}-\ref{Alg: PermutStart-detail2}).

\subsubsection{Upper Bound}\label{sec:UpperBound}

This function calculates an upper bound of the reward for a given node (see Algorithm \ref{Algo:Bounding}). Any node $Q$ in the tree is considered as a leaf node in case $|X^{s_j}_Q|=|Y^{s_j}_Q|=1$ for all $j \in P$. Each leaf node provides a feasible solution $(\hat{x}_{s_1}, \ldots, \hat{x}_{s_p}, \hat{y}_{s_1}, \ldots, \hat{y}_{s_p},\hat{z}_{s_1},\ldots, \hat{z}_{s_p})$ for the problem. In case node $Q$ is not a leaf node (i.e., there exist at least one $j \in P$ such that $|X^j_Q|>1$ or $|Y^j_Q|>1$), lines \ref{Alg:UB:nonleafStart}-\ref{Alg:UB:EndUB} compute the upper bound of the objective value for the subproblem associated with node $Q$. We assume that for each $z \in \Xi=\{\xi_1, \ldots, \xi_m\}$,  $X^z_{ID} := \{x^z_1, \ldots, x^z_{|X^z_{ID}|}\}$ and $Y^z_{ID} := \{y^z_1, \ldots, y^z_{|X^z_{ID}|}\}$ are sorted in the ascending order (Line \ref{Alg:UB:sortCV}). Let $M_z$, $z \in \Xi$, be a matrix whose rows and columns are indexed by $x \in \{x^z_1, \ldots, x^z_{|X^z_{ID}|} \}$ and $y \in \{y^z_1, \ldots, y^z_{|Y^z_{ID}|}\}$. Each entry of this matrix is denoted by $M_z(x,y)$ and it returns the total captured reward for DZs $\cD$ by a SZ with fixed $(x_s,y_s,z_s) = (x,y,z)$. Matrices $M_z$ for all $z \in \Xi$ are created only once while computing initial lower bound using greedy approach and then utilized in each call of \textbf{UpperBound} function. In particular, we store the objective value by solving the PMCLP-PCR with $p=1$ for chosen $z \in \Xi$ using improved PVT algorithm of \cite{BanKia17IJOC}.

To compute an upper bound at node $Q$, we first find the maximum values in $p$ submatrices of $M_z$, $z \in \Xi$, whose rows and columns correspond to $X^{s_j}_Q$ and $Y^{s_j}_Q$, respectively, for $j \in P$ (lines \ref{Alg:UB:StartXandYbar}-\ref{Alg:UB:Ybar}), and then we aggregate these values (lines \ref{Alg:UB:UBupdate1} and \ref{Alg:UB:UBupdate2}). We use matrix $M_z$ for a SZ whose scaling factor $z$ is known at node $Q$, and otherwise, we use $\max_{z \in \Xi}M_z$ for computing the upper bound. Observe that when $X^{s_j}_Q$ contains a single $x$ OSCVs (say $\hat{x}$) such that $z_{s_j} = z \in \Xi$ is known and $\hat{x} \notin X^{z}_{ID}$, we consider a submatrix of $M_z$ with rows and columns corresponding to set $\bar{X} \subset X^{z}_{ID}$ and set $\bar{Y} \subset Y^{z}_{ID}$, respectively. We define the sets $\bar{X}$ and $\bar{Y}$ as follows. In case $x^z_1 < \hat{x} < x^z_{|X^z_{ID}|}$, set $\bar{X}$ is defined by $\{\underline{x},\xb\}$, where $\underline{x}$ and $\xb$ are closest $x$ IDCVs in $X^z_{ID}$ to $\hat{x}$ which are greater and smaller than $\hat{x}$, respectively (line \ref{Alg:UB:SCV1}). Otherwise, set $\bar{X}$ includes either $x^z_1$ (if $\hat{x} < x^z_1$) or $x^z_{|X^z_{ID}|}$(if $\hat{x} > x^z_{|X^z_{ID}|}$) in lines \ref{Alg:UB:SCV2}-\ref{Alg:UB:SCV3}. Likewise, set $\bar{Y} \subset Y^z_{ID}$ is also defined by the same procedure (line \ref{Alg:UB:Ybar}).

\begin{algorithm}[!htbp]
\fontsize{9}{10}\selectfont
 \caption{Upper Bound Calculation}\label{Algo:Bounding}
\begin{algorithmic}[1]
\Function{\fun{UpperBound}}{\pmb{$Q, \cD, \Xi, p$}}
\State $UB \leftarrow 0;$
\label{Alg:UB:UB0}
\State $\{x^z_1,\ldots,x^z_{|X^z_{ID}|} \} \leftarrow$ list of elements in $X_{ID}^z$ sorted in ascending order for any $z \in \Xi$
\label{Alg:UB:sortCV}
  \If {Node $Q$ is a leaf node}{ $UB \leftarrow$ \fun{CoveredReward($Q, \cD, \cS$)}};
  \label{Alg:UB:leafnode}
     \Else
     \label{Alg:UB:nonleafStart}
     			\For{$j \in P$}
     			\label{Alg:UB:StartXandYbar}
     			\State $\hat{z} \leftarrow$ {$\mathcal{Z}(Q,{s_{j}})$};\label{Alg:UB:storeQoS}
     			\If{$\hat{z}=0$} $UB \leftarrow UB+ \max_{z \in \Xi}\big\{M_z(x, y): (x,y) \in X^z_{ID} \times Y^z_{ID}\big\}$
              \label{Alg:UB:UBupdate1}
              \Else 
     			
	           \If{$X_Q^{s_j}$ contains an x OSCV (say $\hat{x}$) and $\hat{x} \notin X^{\hat{z}}_{ID}$}
	           \label{Alg:UB:StartXbar}
	         				\If{$x^{\hat{z}}_1 < \hat{x} < x^{\hat{z}}_{|X^{\hat{z}}_{ID}|}$} $\bar{X} \leftarrow \{\underline{x},\xb\}$ \Comment $\xb$ and $\underline{x}$ are the closest components in $X_{ID}^{\hat{z}}$ to $\hat{x}$ \Statex \hspace{8.5cm} that are larger and smaller than $\hat{x}$, respectively; 	         \label{Alg:UB:SCV1}
	         				\ElsIf{$\hat{x} \leq x^{\hat{z}}_1$}{\hspace{1ex}$\bar{X} \leftarrow \{x^{\hat{z}}_1\}$\;}
               	\label{Alg:UB:SCV2}	\Else{\hspace{1ex}$\bar{X} \leftarrow \{x^{\hat{z}}_{|X^{\hat{z}}_{ID}|}\}$}
               	\label{Alg:UB:SCV3}
	         				\EndIf
	         				
	         		\Else{\hspace{1ex}$\bar{X} \leftarrow X_Q^{s_j}$}
	         		\label{Alg:UB:nonSCV}
	         		\EndIf	
	         		\label{Alg:UB:EndXbar}		
              \State Repeat steps \ref{Alg:UB:StartXbar} to \ref{Alg:UB:EndXbar} for $Y_Q^{s_j}$ to obtain $\bar{Y} \subseteq Y_{ID}$
              \label{Alg:UB:Ybar}
              \State $UB \leftarrow UB+ \max\{M_{\hat{z}}(x, y): (x,y) \in \bar{X} \times \bar{Y}\}$;
              \label{Alg:UB:UBupdate2}
	             
	            \EndIf
	         \EndFor
   \EndIf
   \label{Alg:UB:EndUB}
   \State 	\Return{$UB$}
   \label{Alg:UB:return}
\EndFunction
\end{algorithmic}
\end{algorithm}

For a leaf node, \textbf{UpperBound} calls the \textbf{CoveredReward} function (Algorithm \ref{algo:coveredreward}) to calculate the total reward that has been captured by $p$ SZs with known coordinates and scaling factor. In order to take care of the maximum coverage when multiple SZs with different QoS are overlapping, we need to position SZs with smaller scaling factor (higher QoS) first. Therefore, we sort SZs $\{s_1, \ldots, s_p\}$ based on their scaling factor in the ascending order. Starting with a SZ having smallest scaling factor, the captured reward is calculated by $\sum_{d \in \cD} v^{z_{s_j}}_d A({s_j} \cap d)$ for $j \in P$ and then, by using \textbf{TrimOut} function, we remove the covered parts from the DZs and store a new set of DZs from remaining parts. This procedure is repeated for each fixed SZ and the total calculated covered reward value is aggregated.     
\begin{algorithm}[H]
  \fontsize{9}{11}\selectfont
  \caption{Covered Reward Function} \label{algo:coveredreward}
\begin{algorithmic}[1]
\Function{\fun{CoveredReward}}{\pmb{$Q, \cD, p$}}
    \State $C \leftarrow 0$, $\ccD \leftarrow \cD$;
    \label{algo:cvr: C0}
    \State $\{{s_{i_1}},\ldots, {s_{i_p}}\}$ list of $p$ SZs such that \Statex \hspace{0.5cm} $\mathcal{Z}(Q,{s_{i_1}}) \leq \mathcal{Z}(Q,{s_{i_2}}) \leq \ldots \mathcal{Z}(Q,{s_{i_p}})$;
    \label{algo:cvr:QoSsort}
 	\For{$j \leftarrow i_1$ to $i_p$}
 	\label{algo:cvr:startC}
        \State $(x_{s_j}, y_{s_j}, z_{s_j}) \leftarrow (X^{s_j}_Q, Y^{s_j}_Q, \mathcal{Z}(Q,{s_{j}}))$;
       \For{$d \in \ccD$} $C \leftarrow C + v^{z_{s_j}}_d A({s_j} \cap d)$;
        \State $\ccD \leftarrow TrimOut(\ccD, x_{s_j},y_{s_j},z_{s_j})$;
    \EndFor
    \EndFor
    \label{algo:cvr:EndC}
    \State \Return{$C$}
    \label{algo:cvr:return}
\EndFunction
	\end{algorithmic}
\end{algorithm}

\section{Exact Algorithm for 1D-PMCLP-PC-QoS}\label{sec:1D-PMCLP-PC}

In this section, we provide theoretical properties for the solution space of 1D-PMCLP-PC-QoS where $\Xi_s = \{\xi_s^1\}$ 
for all $s \in \mathcal{S}$, and an exact algorithm to solve it. Recall that in this problem, the DZs (wastewater or trash zones) and SZs (of treatment plants or trash booms) are line segments on $x$-axis (a line representing river). Since scaling factor $z_s$ of each SZ $s$ is fixed but different, i.e., $\xi_s^1$, and $y_{d} = y_s = 0$ for all $d \in \mathcal{D}$ and $s \in \mathcal{S}$, Problem \eqref{eq:MaxRewardFunc} reduces to
\begin{align} \label{eq:MaxReward1D}
\max_{\bx}\bigg\{ f(\bx):=  \sum_{i=1}^n f_i(\bx) = \mathcal{T}_i\lp d_i \cap \lp \cup_{j=1}^p s_j \rp\rp \bigg\},
\end{align}
where $\mathcal{T}_i(.)$ returns the total reward from line segment captured by $s_1, \ldots, s_p$. In other words, $\bz = \{\xi^1_{s_1}, \ldots, \xi^1_{s_p}\}$ and $\by = {\bf 0}$ in \eqref{eq:MaxRewardFunc}. Assume that the base SZ $s_0$ for 1D-PMCLP-PC-QoS is a line segment with known width $w_{s_0}$. Note that the reward rates for covering each DZ $d_i$ are reward per unit length. For $p=1$, 1D-PMCLP-PC-QoS is a special case of the PMCLP-PCR and therefore, the properties for solution space and exact algorithms provided by \cite{SonStaGol06} and \cite{BanKia17IJOC} for PMCLP-PCR are directly applicable to 1D-PMCLP-PC. However, the same is not true for $p\geq 2$, primarily because each SZ has different width. 

\subsection{Solution Space of 1D-PMCLP-PC-QoS}

We establish relation between PMCLP-PCR-QoS with $\Xi_s = \{\xi_1, \ldots, \xi_m\}$ for all $s \in \mathcal{S}$ and 1D-PMCLP-PC-QoS with $\Xi_s = \{\xi_s^1\}$ for all $s \in \mathcal{S}$. Notice that in case $\xi_s^1 \in \{\xi_1, \ldots, \xi_m\}$ for all $s \in \mathcal{S}$, then the objective function of 1D-PMCLP-PC-QoS is conceptually same as the objective function of PMCLP-PCR-QoS with fixed $(\by = {\bf 0},\bz)$ for SZs in the 1D solution space. In other words, the objective function of 1D-PMCLP-PC-QoS is a piecewise linear function with break points at $x$ DCVs and SCVs. Since each SZ has different width in 1D-PMCLP-PC-QoS, we define DCVs associated to each SZ $s_j$, $j\in \{1,\ldots,p\}$, by $X^{z_{s_j}}_D = \cup_i X^{d_i, z_{s_j}}_D$ (see Definition \ref{def:DCVs})
and denote the set of all IDCVs and ODCVs associated to SZ $s_j$ by $X^{z_{s_j}}_{ID}$ and $X^{z_{s_j}}_{OD}$, respectively. Each SZ with a fixed position creates a set of three or four SCVs for other SZs. The set of SCVs generated by SZ $s_j$, positioned at $x_{s_j}$, for another SZ $s_k$ ($k\neq j$) is defined by $\mathcal{X}^{s_j, z_{s_k}}_S$ (see Definition \ref{def:SCVs}). Likewise, we define $\mathcal{X}^{z_s}_S(\mathcal{J})$, $\mathcal{X}^{z_s}_{IS}(\mathcal{J})$, and $\mathcal{X}^{z_s}_{OS}(\mathcal{J})$ for $s \in \mathcal{S}\bs\mathcal{J}$. Since the scaling factor of all SZs is fixed, there is no QCV, i.e., $X^{d_i}_Q = \emptyset$ for all $i$. Furthermore, Observations \ref{obs:ODCVnISCV} and \ref{obs:gxbreakpoints} are valid for the objective function of 1D-PMCLP-PC-QoS.

\begin{theorem}\label{theorem:1D}
There exists an optimal solution $\bx^*=(x_{s_1}^*, \ldots,x_{s_p}^*) \in \sr^p$ of the 1D-PMCLP-PC-QoS 
and a non-repetitive sequence $s_{t_1}, s_{t_2}, \ldots, s_{t_p}$ of the SZs such that $x_{s_{t_1}}^* \in X^{z_{s_{t_1}}}_{ID}$ and for $k=2,\ldots,p$, $x_{s_{t_k}}^* \in X^z_{ID} \cup \mathcal{X}^{z}_{OS}(T_k)$ where $z = \xi^1_{s_{t_k}}$, $T_k=\{t_1, \ldots, t_{k-1}\}$, and $t_k \in P$.
\end{theorem}

\noindent {\it Proof.} This theorem can be proved using the same arguments as used in the proof of Theorem~\ref{theorem:MSP-MQoS} for condition ($i$), except that $z^*_s =\hat{z}_s = \xi^1_s$ for all $s \in S$ and $X_Q = \emptyset$ for all $i$.

\subsection{Exact Algorithm for 1D-PMCLP-PC-QoS}
The 1D-PMCLP-PC-QoS for $p=1$ and $p \geq 2$ with $\Xi_{s} = \{1\}$ for all $s \in \mathcal{S}$, are special cases of the PMCLP-PCR, and therefore, the algorithms provided by \cite{SonStaGol06} and \cite{BanKia17IJOC}, respectively, are directly applicable for these special cases. However, when $p\geq 2$ and each SZ has different scaling factor, these approaches do not guarantee to provide an optimal solution for 1D-PMCLP-PC-QoS. Assuming that $\xi_s^1 \in \Xi = \{\xi_1, \ldots, \xi_m\}$ for all $s \in \mathcal{S}$, our proposed exact algorithm for PMCLP-PCR-QoS with $\by = {\bf 0}$, $\Xi:=\{\xi^1_{s_1}, \ldots,\xi^1_{s_p}\}$, and $z_{s_j} = \xi^1_{s_j}$, $j \in P$, can be applied to find an optimal solution for 1D-PMCLP-PC-QoS. In this section, we provide another more computationally efficient B\&B based exact algorithm for the 1D-PMCLP-PC-QoS. The main body and upper bound computation in this algorithm are similar to Algorithms \ref{Algo:MainBody} and \ref{Algo:Bounding} with some minor modifications (discussed below), but the branching routine involves features pertinent to this problem. 

\subsubsection{Main Body and Upper Bound} The outline of this approach is also given by Algorithm \ref{Algo:MainBody}, where node $Q$ stores $p$ subsets of $x$ CVs (IDCVs and OSCVs), $X^{s_j}_{Q}$, $j \in P$, and an indicator $\mathcal{BS}(Q)$. Each node represents a restricted 1D-PMCLP-PC-QoS problem where SZ $s_j \in \mathcal{S}$ has $x_{s_j} \in X^{s_j}_Q$. In this B\&B tree, a leaf node $Q$ has $|X^{s_j}_Q|=1$ for all $j \in P$, and at root node $Q_0$, set $X^{s}_{Q_0} = X^{\xi^1_s}_{ID}$ for all $s \in \mathcal{S}$. In order to compute an upper bound ($\UB$) on the objective function of the restricted subproblem at node $Q$, we utilize and modify Algorithm \ref{Algo:Bounding} as follows. For a leaf node, we call \textbf{CoveredReward} function (Algorithm \ref{algo:coveredreward}) that provides a feasible solution $(\hat{x}_{s_1}, \ldots, \hat{x}_{s_p})$ for this problem as well. In case node $Q$ is not a leaf node, an upper bound is computed using Lines \ref{Alg:UB:nonleafStart}-\ref{Alg:UB:EndUB} of Algorithm \ref{Algo:Bounding}, where $Y^z_{ID} = \bar{Y} = \{0\}$ for all $z$, $\Xi = \{\xi^1_{s_1}, \ldots, \xi^1_{s_p}\}$, and $\mathcal{Z}(Q, s_j) = \xi^1_{s_j}$, $j \in P$. As a result, $M_z$, $z \in \Xi$, becomes a row vector whose rows are indexed by $x \in \{x^z_1, \ldots, x^z_{|X^z_{ID}|} \}$ and each entry of this vector is denoted by $M_z(x,0)$ that returns the total captured reward for DZs $\cD$ by a SZ with fixed $(x_s,y_s,z_s) = (x,0,z)$. 

\subsubsection{Branching} This function generates a list of child nodes for a given parent node $Q$. These child nodes are added to the \texttt{LCN} and the most left child node of $Q$ is selected from \texttt{LCN} in the next iteration. Since the scaling factor of each SZ comes from non-homogeneous sets of scaling factors ($\Xi_{s_t} \neq \Xi_{s_w}$ for $t \neq w$), then according to Theorem \ref{theorem:1D} we have to consider all sequences of SZs, $t_1, t_2, \ldots, t_p \in P$ in the branching. Observe that the symmetries in the B\&B tree for PMCLP-PCR-QoS that resulted in the reduction of sequences of branching SZs along an axis, no longer exist for 1D-PMCLP-PC-QoS. We introduce a different strategy to truncate B\&B tree for the latter, thereby leading to a new branching routine (Algorithm \ref{Algo:Branch1D}), denoted by \fun{Branching1D}($Q$). For root node $Q=Q_0$, we create $p$ child nodes (Lines \ref{line:rootnode1}-\ref{line:rootnode2} of Algorithm \ref{Algo:Branch1D}) where each child node $T$ has same properties as $Q_0$, except that indicator $\mathcal{BS}(T) \in \mathcal{S}$ is different for all child nodes. We refer to these nodes as first-level branching nodes, denoted by $\overline{Q}^j_0$, $j \in P$, where $index(\mathcal{BS}(\overline{Q}^j_0)) = j$. For each node $Q$ in the subtree with $\overline{Q}^j_0$, $j \in P$, as root node, we define another indicator $\texttt{BSFL}(Q)$ to store $index(\mathcal{BS}(\overline{Q}^j_0))$, i.e., $j$, and harness this information to reduce the number of nodes in the B\&B tree.

\begin{algorithm}[!h]
\fontsize{9}{9}\selectfont
\caption{Branching for 1D-PMCLP-PC-QoS}
\label{Algo:Branch1D}
\begin{algorithmic}[1]
\Function{\fun{Branching1D}}{\pmb{$Q$}}
\If{$Q$ is the root node $Q_0$} \label{line:rootnode1}
\For{$j \in P = \{1,\ldots,p\}$}
\State $T\leftarrow$ {$NewChildNode(Q)$}
\State $\mathcal{BS}(T)\leftarrow$ {$s_j$}
\State $\texttt{BSFL}(T)\leftarrow$ {$j$}
\EndFor\label{line:rootnode2}
\Else
\If{$|X_Q^{\mathcal{BS}(Q)}|=1$}
\State $L \leftarrow \{k \in P: |X^{s_k}_Q| =1\}$ \label{line:revisedpermutation1}
\For{$l \in P \setminus L$}
\For{$\hat{x} \in \X^{s_l}_{OS}(L)$}
\State $T\leftarrow$ {$NewChildNode(Q)$}
\State $\mathcal{BS}(T)\leftarrow$ {$s_l$} 
\State $X_T^{\mathcal{BS}(T)} \leftarrow \hat{x}$ 
\EndFor
\If{$l > \texttt{BSFL}(Q)$}
\label{Alg: 1DB: PermutStart}
\State $T\leftarrow$ {$NewChildNode(Q)$}
\State $\mathcal{BS}(T)\leftarrow$ {$s_l$} 
\State $X_T^{\mathcal{BS}(T)} \leftarrow X_{ID}^{\xi_{s_l}^1}$
\EndIf \label{line:revisedpermutation2}
\EndFor
\Else
\For{$n_x \leftarrow 1$ to $NumPartitions(X^{\mathcal{BS}(Q)}_Q)$} \label{line:partition1}
\State $T\leftarrow$ {$NewChildNode(Q)$}
\State $X_T^{\mathcal{BS}(T)} \leftarrow$ {$GetPartition(X_Q^{\mathcal{BS}(Q)},n_x)$} \label{line:partition2}
\EndFor
\EndIf
\EndIf
\EndFunction
	\end{algorithmic}
\end{algorithm}

In case $|X^{\mathcal{BS}(Q)}_Q| \neq 1$ at a node $Q$, function \fun{Branching1D}($Q$) performs the partition-based branching  in Lines \ref{line:partition1}-\ref{line:partition2} (also discussed in Section \ref{sec:Branching}). However, whenever $|X^{\mathcal{BS}(Q)}_Q| = 1$ at a node $Q$, i.e., SZ $\mathcal{BS}(Q)$ has fixed $x$ coordinate, function \fun{Branching1D}($Q$) conducts a revised permutation branching (Lines \ref{line:revisedpermutation1}-\ref{line:revisedpermutation2}) that works as follows. For 1D-PMCLP-PC-QoS with $p=2$, according to Theorem \ref{theorem:1D}, an optimal solution $(x^*_{s_1}, x^*_{s_2})$ belongs to the set {\fontsize{11}{11}\selectfont
\begin{align}\label{eq:set1DPMCLPp=2}
\bigg\{(x_{s_1}, x_{s_2}): x_{s_{t_1}} \in X^{z_{s_{t_1}}}_{ID} \text{ and } x_{s_{t_2}} \in X^{z_{s_{t_2}}}_{ID} \cup \mathcal{X}^{z_{s_{t_2}}}_{OS}(\{t_1\}) \text{ for all } (t_1, t_2) \in \{(1,2), (2,1)\}\bigg\}.
\end{align}}
We observe that 
$\big\{\big(x_{s_{t_1}}, x_{s_{t_2}}\big) \in X^{z_{s_{t_1}}}_{ID} \times X^{z_{s_{t_2}}}_{ID}\big\}$ is same for each pair $(t_1, t_2)\in \{(1,2),(2,1)\}$ and hence, considering them once reduces the number of nodes in the B\&B tree. Based on this observation, set \eqref{eq:set1DPMCLPp=2} that can be rewritten as
\begin{align*}
\bigg\{(x_{s_1}, x_{s_2}): X^{z_{s_1}}_{ID}  \times \lp X^{z_{s_2}}_{ID} \cup \mathcal{X}^{z_{s_2}}_{OS}(\{1\}) \rp \bigg\} \bigcup \bigg\{(x_{s_2}, x_{s_1}): X^{z_{s_2}}_{ID}  \times \lp X^{z_{s_1}}_{ID} \cup \mathcal{X}^{z_{s_1}}_{OS}(\{2\}) \rp \bigg\} \end{align*}
 is equivalent to
\begin{align*}
\bigg\{(x_{s_1}, x_{s_2}): X^{z_{s_1}}_{ID}  \times \lp X^{z_{s_2}}_{ID} \cup \mathcal{X}^{z_{s_2}}_{OS}(\{1\}) \rp \bigg\} \bigcup \bigg\{(x_{s_2}, x_{s_1}): X^{z_{s_2}}_{ID}  \times \mathcal{X}^{z_{s_1}}_{OS}(\{2\}) \bigg\}.
\end{align*}
We incorporate this observation in revised permutation branching by considering: (a) $x_{s_1} \in X^{z_{s_1}}_{ID}$ and $x_{s_2} \in X^{z_{s_2}}_{ID} \cup \mathcal{X}^{z_{s_2}}_{OS}(\{1\})$ in subtree with $\overline{Q}^{1}_0$ as root node, and (b) $x_{s_2} \in X^{z_{s_2}}_{ID}$ and $x_{s_1} \in \mathcal{X}^{z_{s_1}}_{OS}(\{2\})$ in subtree with $\overline{Q}^{2}_0$ as root node. In general, for each subtree with $\overline{Q}^{j}_0$, $j \in P$, as root node we set $s_{t_1} = index(\mathcal{BS}(\overline{Q}^{j}_0)) = j$ and explore the solution space 
\begin{align*}
    \bigg\{\lp x_{s_{t_1}}, x_{s_{t_2}}, \ldots, x_{s_{t_p}} \rp:  x_{s_{t_1}} \in X^{z_{s_{t_1}}}_{ID}, \ \ & x_{s_{t_k}} \in \mathcal{X}^{z_{s_{t_k}}}_{OS}(\{T_k\}) \text{ if } t_k < t_1,  \\ &   x_{s_{t_k}} \in X^{z_{s_{t_k}}}_{ID} \cup \mathcal{X}^{z_{s_{t_k}}}_{OS}(\{T_k\}) \text{ if } t_k > t_1\bigg\}
\end{align*}
for all non-repetitive sequences of $t_2, \ldots, t_p \in P \bs \{t_1\}$. Accordingly, we develop the B\&B tree using the revised permutation branching (Lines \ref{line:revisedpermutation1}-\ref{line:revisedpermutation2}). More specifically, whenever $|X^{\mathcal{BS}(Q)}_Q| = 1$ at a node $Q$, we define a set $L$ of SZ indices whose $x$ coordinate is fixed and for each $l \in P\bs L$, we create $|\X^{s_l}_{OS}(L)|$ number of child nodes using the \funi{NewChildNode} routine where each child node $T$ has $\mathcal{BS}(T) = l$ and singleton $X^{\mathcal{BS}(T)}_T$ contains an $x$ OSCV $\hat{x}\in \X^{s_l}_{OS}(L)$. In addition, for each $l \in P \bs L$ such that $l > \texttt{BSFL}(Q)$, we also create a child node $T$ of node $Q$ that has $\mathcal{BS}(T) = l$ and $X^{\mathcal{BS}(T)}_T$ contains $x$ IDCVs corresponding to SZ $s_l$, i.e., $X^{\xi^1_{s_l}}_{ID}$. As mentioned before, this technique reduces in the number of nodes of the B\&B tree, thereby leading to a computationally efficient solution approach for 1D-PMCLP-PC-QoS.

\section{Computational Experiments}\label{sec:ComExp}

For computational study, we perform three sets of experiments to evaluate the effectiveness and efficiency of the solution approaches presented in this paper. Specifically, in the first set of the experiments, we compare the performance of our proposed exact method for the PMCLP-PCR (i.e., PMCLP-PCR-QoS with fixed and same scaling factor $\Xi = \{1\}$), with the existing algorithm proposed by \cite{BanKia17IJOC}. The second and third sets of experiments are conducted to evaluate the performance of our proposed exact and approximation algorithms for PMCLP-PCR-QoS with $\Xi = \{\xi_k = k\}_{k=1}^m$ and for 1D-PMCLP-PC-QoS with $\Xi_{s_j} = \{\xi^1_{s_j} = j\}$ for all $j \in P$. We assume that $\eta(z) = z$. All three sets of experiments are conducted by solving randomly generated instances. We implemented these experiments in the Python 3.5.5 and ran them on a Dell Precision 5820 workstation with 3.80 GHz Intel Core i7-9800 Processor with 32.0 GB RAM and Windows 10.

\subsection{Instance Generation Procedure}

To generate instances for our computational experiments, we adopt the instance generation procedure of \cite{SonStaGol06,BanKia17IJOC}. We consider a square region of size $1000 \times 1000$ in which the DZs are located. Each rectangular DZ, $d \in \cD$, is specified by lower left corner coordinates $(x_d, y_d)$, width $w_d$, length $l_d$, and reward rate $v_d$. The coordinates of lower left corner ($x_d, y_d$) of each DZ is generated as follows. First, we randomly generate three center points in the mentioned square region using a uniform distribution. These center points are referred to as concentration points and we associate a circular region of radius $r = 270$ around each center point. To randomly generate and position each DZ, we first decide whether it is "anchored" to a concentration point (with probability of 0.31 for each concentration point), or it is "free" (with probability of 0.07). The lower left corner of an anchored DZ is randomly located within the circular region assigned to the corresponding concentration point, while the free DZ is randomly located anywhere in the entire square region ($1000 \times 1000$). This DZ positioning approach is used to mimic the real-word situation where the large percentage of demand is concentrated around some population centers, while a small percentage is scattered all over the region. If the DZ is anchored, the $x$ (and $y$) coordinate of its lower left corner are generated from a normal distribution, where mean is the $x$ (and $y$) coordinate of the corresponding concentration point and standard deviation is $r/3$. If the DZ is free, it will be placed in the square region using a uniform distribution. The width and length of each DZ follow a uniform distribution with $w_d, l_d \sim \text{uniform}[5,50]$. The reward rate is also generated based on a uniform distribution where $v_d \sim \text{uniform}[1,10]$. Also, we consider the base SZ $s_0$  of dimensions $(w_{s_0},h_{s_0})= (50,40)$. The instance generation for 1D-PMCLP-PC-QoS follows the same procedure on a line segment of length $1000$ units on $x$-axis where $w_d = 0$ and $y_d = 0$ for all $d \in \cD$, and $w_{s_0}=0$. 
%


\subsection{Computational Results for PMCLP-PCR}
In Table \ref{tab:PMCLP-PCR}, we present results of our first set of computational experiments where we compare the performance of our new exact algorithm (that incorporates Theorem \ref{theorem:MSP-MQoS}) with an existing (benchmark) algorithm of \cite{BanKia17IJOC} for PMCLP-PCR.
Each row in Table \ref{tab:PMCLP-PCR} represents the average over 10 instances. In this table, we use $N$, $T$, and $T_1$ to denote the total number of branch-and-bound nodes, total solution time (in seconds), and time taken in seconds to find an optimal solution, respectively. The subscripts $\mathcal{S}$ and $\mathcal{OS}$ denote the results for the existing algorithm \citep{BanKia17IJOC} and our proposed algorithm, respectively. These subscripts represent that in the former all the SCVs are considered in the solution search space, whereas in the latter the solution space is reduced to only OSCVs. Table \ref{tab:PMCLP-PCR} shows that our proposed method always has lesser number of nodes traversed, time taken to solve, and time taken to find an optimal solution, compared to the benchmark. The last three columns provides percentage improvement in $N$, $T$, and $T_1$ that are quantified by $\texttt{ImprN}:=100\times\lp 1 - \frac{N_{\mathcal{OS}}}{N_{\mathcal{S}}}\rp$, $\texttt{ImprT}:=100\times\lp 1 - \frac{T_{\mathcal{OS}}}{T_{\mathcal{S}}}\rp$, and $\texttt{ImprT}_1:=100\times \lp 1 - \frac{T_{1\mathcal{OS}}}{T_{1\mathcal{S}}}\rp$, respectively. 
Notice that for $p=2$, $3$, and $4$, we have in average $32.06\%$, $33.23\%$, and $41.36\%$ improvements, respectively, in terms of number of traversed nodes ($N$) as well as $35.67\%$, $26.89\%$, and $39.67\%$ average improvements, respectively, in terms of total time ($T)$. Both $N$ and $T$ have maximum reduction of $67\%$ and $64.9\%$, respectively, using our proposed method. However, the improvement in terms of the time taken to find an optimal solution ($T_1$) has more variation among different number of DZs (varies between $0.0\%$ and $87.5\%$). For instance, when $p=2$, the $T_1$ improves in average $35.4\%$, and when $p=3$ and $p=4$, we get even better average reduction of $58.25\%$ and $51.43\%$, respectively.

\begin{table}[t]\footnotesize
\centering 
\setlength{\tabcolsep}{0.75pt}\renewcommand{\arraystretch}{1.4}
\addtolength{\tabcolsep}{-3.6pt}
\centering
\caption{{Computational results for PMCLP-PCR }} \label{tab:PMCLP-PCR}\small
\begin{tabular}{c|c|c|c|c|c|c|c|c|c|c}
\hline

\multirow{2}{*}{  \ $p$ \   } &
  \multirow{2}{*}{ \ $n$ \  } &
  \multicolumn{3}{c|}{This Paper} &
  \multicolumn{3}{c|}{\cite{BanKia17IJOC}} &
  \multicolumn{3}{c}{Improvement (\%)}
   \\ \cline{3-11}
 &
   &
  \multirow{2}{*}{$N_{\mathcal{OS}}$} &
  \multirow{2}{*}{\ \  $T_{\mathcal{OS}}$(s) \ }&
  \multirow{2}{*}{\ \  $T_{1_{\mathcal{OS}}}$(s) \ }&
  \multirow{2}{*}{$N_{\mathcal{S}}$} &
  \multirow{2}{*}{ \ \ $T_{\mathcal{S}}$(s) \ } &
  \multirow{2}{*}{ \ \ $T_{1{\mathcal{S}}}$(s) \ }&
  \multirow{2}{*}{\ \ \texttt{ImprN} \ } & \multirow{2}{*}{ \ \ \texttt{ImprT } } &
  \multirow{2}{*}{\ \ $\texttt{ImprT}_1$} \\

&  &  &  &  & &  & &  &  & \\
\hline
\multirow{5}{*}{2} 
& 10 & 29 & 0.013 & 0.000 & 88 & 0.037 & 0.000  & 67.0 & 64.9 & 0.0\\
& 20 & 201 & 0.095 &0.010 & 303 & 0.159 &0.044  & 33.6 & 40.0 & 77.2\\
& 30 & 297 & 0.187 & 0.055 & 356 & 0.226 & 0.082 & 16.6 & 17.3 & 32.9\\
& 50 & 564 & 0.449 & 0.075 & 666 & 0.604 & 0.098 & 15.3 & 25.7 & 48.9\\
& 70 & 651 & 0.674 & 0.094 & 976 & 1.028 & 0.105 & 33.3 & 34.5 & 10.4\\
& 100 & 3,601 & 5.448 & 1.253 & 4,843 & 7.962 & 1.31 & 26.6 & 31.6 & 4.3\\
\hline
\multirow{5}{*}{3} 
& 10 & 2,801 & 1.371 & 0.021 & 3,806 & 2.053 & 0.023  & 26.4 & 33.3 & 9.5\\
& 20 & 9,395 & 7.457 & 0.753 & 17,748 & 13.98 & 1.427 & 47.0 & 46.7 & 47.2\\
& 30 & 62,563 & 66.61 & 3.92 & 86,697 & 85.60 & 9.37 & 27.8 & 22.1 & 58.1\\
& 50 & 109,065 & 155.6 & 9.024 & 146,582 & 216.6 & 32.76 & 25.6 & 28.4 & 72.4 \\
& 70 & 237,434 & 354.6 & 5.742 & 290,098 & 389.3 & 28.73 & 18.1 & 8.99 & 80.0\\
& \ \ 100 \ \ & 545,852 & 1,311 & 31.54 & \ $1.21 \times 10^6$ \  & 2,731 & 178.7 & 54.5 & 51.9 & 82.3\\
\hline
\multirow{3}{*}{4} 
& 10 & 894,255 & 720.4 & 0.012 & $1.9 \times 10^6$ & 1,782  & 0.096 & 53.9 & 59.5 & 87.5\\
& 20 &  \ \ $1.4 \times 10^6$ \  \ &  1,904 &  \ 397.3 \ & $2.2 \times 10^6$ &  \ 2,513 \ & \ 441.5 \ & 34.9 & 24.2 & 10.0\\
& 30 & $6.9 \times 10^6$ & \ \ 10,390 \ \ & \ \ \ 1,002 \ \ \ &  \ \ $10.7 \times 10^6$ \ \ &  \ \ \ 16,132 \ \ \ & \ \ \ 2,321 \ \ \ & 35.3 & 35.6 & 56.8\\
\hline
\end{tabular}

\end{table}

\subsection{Computational Results for PMCLP-PCR-QoS}

Through the results of the second set of experiments, we analyze the performance of our proposed greedy approximation and branch-and-bound exact algorithms for PMCLP-PCR-QoS with $|\Xi|=m \geq 2$. These experiments are conducted on randomly generated instances with number of SZs $p \in \{2,3,4\}$. The results are presented in Table \ref{tab:PMCLP-PCR-QoS}, where each row provides an average over results for 10 instances. Columns labelled as $T$ and $T_1$ report the total solution time and time taken to find an optimal solution, respectively. Both $T$ and $T_1$ are reported in seconds. In this experiment,  we set a time limit of 5 hours (18,000 seconds). Therefore, the algorithm terminates after 5 hours and reports the best solution found in these 5 hours. If we cannot guarantee that we have found an optimal solution within the mentioned time limit, we will use $18,000+$ in the $T$ column and '-' in the $T_1$, $T_1/T$, and $\alpha$ columns. Also, we report the total number of nodes traversed in our proposed B\&B tree at the termination.

\begin{table}[!htbp]\footnotesize
\centering 
\setlength{\tabcolsep}{8pt}\renewcommand{\arraystretch}{1.4}
\caption{Computational results for PMCLP-PCR-QoS with $m \geq 2$}\label{tab:PMCLP-PCR-QoS}

\begin{tabular}{c|c|c| c| c|c |c |c|c}
\hline
$p$ & $m$ & \# {of} DZs &  \# {of} Nodes  & $T (s)$ & $T_1 (s)$ & $T_1/T$  & $T_H (s)$&$\alpha$ \\ 
\hline
\multirow{9}{*}{2} &

\multirow{3}{*}{2} & 10 & 96  & 0.040 & 0.000 & 0.000 & 0.042 & 0.999\\ &  & 50 & 2,915 & 3.042 & 0.250 & 0.082 & 3.251 & 0.996\\ & & 100 & 13,092 & 19.19 & 5.681 & 0.295 & 23.82& 0.990 \\

\cline{2-9}
 & \multirow{3}{*}{3} & 10 & 223  & 0.083 & 0.000 & 0.000 & 0.072 & 0.999 \\ &  & 50 & 7,983 & 7.082 & 1.436 & 0.202 & 5.712 & 0.986\\ &  & 100 & 23,379 & 37.63 & 5.347 & 0.142 & 35.45 & 0.993\\

\cline{2-9}
 & \multirow{3}{*}{4} & 10 & 857  & 0.346 & 0.018 &  0.052  & 0.124 & 0.993\\ &  & 50 & 12,292 & 12.70  & 1.934 & 0.152 & 6.138 & 0.990\\ &  & 100 & 37,655 & 69.37 & 12.89 & 0.185 & 47.81 & 0.984\\
 
\cline{2-9}
 & \multirow{3}{*}{5} & 10 & 2,121  & 0.803 & 0.137 &  0.170  & 0.13 & 0.996\\ &  & 50 & 37,019 & 33.27 & 2.94 & 0.088 & 7.71 & 0.999\\ &  & 100 & 66,661 & 113.64 & 34.59 & 0.268 & 57.64 & 0.981\\
\hline

\multirow{9}{*}{3} &

\multirow{3}{*}{2} & 10 & 51,803 & 32.12 & 4.27 & 0.132 & 0.063 & 0.992\\ &  & 25 & 215,453 & 210.42 & 33.93 & 0.161 & 7.322 & 0.982\\&  & 50 & 326,283 & 422.35 & 129.83 & 0.307 & 49.81 & 0.974\\

\cline{2-9}
 & \multirow{3}{*}{3} & 10 & 77,092 & 38.84 & 0.953 & 0.024 & 0.122 & 0.989\\ &  & 25 & 511,442 & 442.15 & 103.26 & 0.233 & 11.28 & 0.967\\ &  & 50 & $2.2\times 10^6$ & 3,342.2 & 403.29 & 0.120 & 74.351 & 0.985\\
\cline{2-9}
 & \multirow{3}{*}{4} & 10 & 550,522 & 287.60 & 43.54 & 0.151 & 0.127 & 0.995\\ &  & 25 & $2.2 \times 10^6$ & 1,968.9 & 193.50 & 0.098 & 13.02  & 0.996\\ &  & 50 & $5.8\times 10^6$  & 9,248.8 & 1,679 & 0.181 & 94.64 & 0.973\\
 
\cline{2-9}
 & \multirow{3}{*}{5} & 10 & $1.9\times 10^6$ & 1,161.8 & 258.16 & 0.222 & 0.169  & 0.979\\ &  & 25 & $4.2 \times 10^6$ & 3,430.3 & 429.84 & 0.125 & 15.08 & 0.987\\ &  & 50 & $8.4 \times 10^6$  & 12,476 & 1965.9 & 0.157 & 114.5 & 0.986\\
\hline

\multirow{9}{*}{4} &

\multirow{3}{*}{2} & 10 & $4.1 \times 10^6$ & 3,287.1 & 196.05 & 0.059 & 0.137 & 0.988 \\ &  & 25 & $9.6 \times 10^6$ & 11,827 & 350.70 & 0.029 & 10.95 & 0.995\\ &  & 50 & $10.2 \times 10^6$ & 18,000+ & -- & -- & 66.41 & --\\
\cline{2-9}
 & \multirow{3}{*}{3} & 10 & $6.9 \times 10^6$  & 4,818.7&0.288 & 0.058 & 0.161 & 0.987 \\ &  & 25 & $12.8 \times 10^6$  & 14,696& 1,218.4 & 0.082 & 12.13 & 0.982\\ &  & 50 & $19.3 \times 10^6$ & 18,000+ & -- & -- & 98.91 & --\\

\cline{2-9}
 & \multirow{3}{*}{4} & 10 & $9.9 \times 10^6$ & 14,570 & 1,321.3 &  0.09 & 0.184 & 0.986\\ &  & 25 & $13.7 \times 10^6$ & 15,085  & 1,776.7 & 0.118 & 17.64 & 0.973\\ &  & 50 & $21.8 \times 10^6$ & 18,000+ & -- & -- & 121.7 & --\\
 
\cline{2-9}
 & \multirow{3}{*}{5} & 10 & $14.8 \times 10^6$ & 16,873 & 569.25 &  0.033  & 0.221 & 0.998\\ &  & 25 & $22.5 \times 10^6$ & 18,000$+$ & -- & -- & 21.63 & --\\ &  & 50 & $31.2 \times 10^6$  & 18,000$+$ & -- & -- & 145.6 & --\\
 
\hline

\end{tabular}
\end{table}

Table \ref{tab:PMCLP-PCR-QoS} illustrates that the difficulty of PMCLP-PCR-QoS substantially increases as the number of SZs ($p$) increases. We have already shown that PMCLP-PCR-QoS is NP-hard when $p$ is a part of the input. However, our proposed method can solve relatively large instances in a short amount of time. For example, for $p=2$, $n=10$ (number of DZs), and $m\in\{2,3,4,5\}$, the total solution time ($T$) of our proposed algorithm is less than $1$ seconds. Also, for $p=2$, $n \leq 100$, and $m\in\{2,3,4\}$, the total time is less than $1$ minutes. When $p=3$, $n \leq 50$, and $m=2$, the maximum total time is around $7$ minutes. Moreover, the time taken by our proposed method to solve the provided instances is greatly smaller than the time taken by the explicit enumeration of all CVs. For example, in case $p =2$, $m=2$, and $n=100$, we need to traverse $1600 \times 1600 = 2,560,000$ number of nodes in the explicit enumeration, whereas our proposed method find an optimal solution by only traversing $13,092$ number of nodes in the B\&B tree. In other words, our proposed algorithm benefits from the implicit enumeration of breakpoints by concentrating on areas with high rewards, and using good lower bounds that are provided by our proposed greedy algorithm. The column labelled as $T_1/T$ provides the ratio of $T_1$ and $T$. Based on Table \ref{tab:PMCLP-PCR-QoS}, the values of the column $T_1/T$ are quite small (i.e., mostly less than 0.2). This means that because of starting from a good lower bound, our proposed algorithm reaches an optimal solution fast and spends the remaining time to prove the optimality of this solution.

One of the advantages of proposing the exact algorithm for PMCLP-PCR-QoS is to evaluate the performance of any heuristic/approximation algorithm for it. In Table \ref{tab:PMCLP-PCR-QoS}, the column $T_H$ refers to the total time taken by the proposed greedy algorithm to solve PMCLP-PCR-QoS instances. As we see, the $T_H$ is smaller than the $T$ as number of SZs ($p$) and number of possible scaling factors ($m$) increase. The column labelled as $\alpha$ provides the empirical approximation ratio of the reward captured by the greedy algorithm to the optimal covered reward by the exact algorithm. Based on Theorem \ref{thm:GreedyAlgo}, the approximation ratio of the proposed greedy solution is $1 -\frac{1}{e}$ (almost $63.2\%$). However, after finding an optimal solution, we can ensure that the greedy algorithm suboptimal solution has a comparably good quality. For example, for the instances considered in Table \ref{tab:PMCLP-PCR-QoS}, this ratio is at least $96.9\%$. In conclusion, Table \ref{tab:PMCLP-PCR-QoS} will help in deciding whether to use the greedy or the exact algorithm depending on the input parameters of PMCLP-PCR-QoS, and the desired tolerance and available computational resources. For example, when $p=2$, the difference between $T_H$ and $T$ is few seconds, as a result, we prefer to use the exact algorithm and reach to the optimal coverage. However, for $p=4$, the gap between $T_H$ and $T$ is huge, thus, one may prefer to use faster greedy algorithm and find a good suboptimal coverage instead of an optimal solution. Instances for which we could not obtain optimal solution, we cannot guarantee the high quality of the greedy solution because we cannot compute $\alpha$.

\subsection{Computational Results for 1D-PMCLP-PCR-QoS}

In the third set of experiments, we plan to evaluate the performance of our proposed greedy approximation and branch-and-bound exact algorithms for 1D-PMCLP-PC-QoS with $\Xi_s = \{\xi_s\}$, for all $s \in \mathcal{S}$. These experiments are performed on the randomly generated instances with number of SZs $p \in \{2,3,4\}$, and the results of the experiments are reported in Table \ref{tab:1D-PMCLP-PC-QoS}. In Table \ref{tab:1D-PMCLP-PC-QoS}, the columns labelled as $T$, $T_1$, $T_1/T$, $T_H$, and $\alpha$ are defined similar to the Table \ref{tab:PMCLP-PCR-QoS}. When an optimal solution cannot be found within time limit of 18000 seconds, we use $18,000+$ in the $T$ column and `-' in the $T_1$, $T_1/T$, and $\alpha$ columns. Based on Table \ref{tab:1D-PMCLP-PC-QoS}, our proposed exact method can solve the 1D-PMCLP-PC-QoS instances in less than $1$ second when $p=2$. Also, in case $p=3$, our algorithm solved the instances with number of DZs $n \in\{10,20,30,50,70\}$ in less than $1$ minutes. However, when $p$ increases (e.g., $p=4$) the total solution time ($T$) grows exponentially.

\begin{table}[!htbp]\footnotesize
\centering 
\setlength{\tabcolsep}{6pt}\renewcommand{\arraystretch}{1.4}
\caption{Computational results for 1D-PMCLP-PCR-QoS}\label{tab:1D-PMCLP-PC-QoS}
\centering
\begin{tabular}{ c|c|c| c|c |c |c|c}
\hline
$p$ &  \# {of} DZs &  \# {of} Nodes  & $T (s)$ & $T_1 (s)$ & $T_1/T$  & $T_H (s)$ & $\alpha$ \\ 
\hline
\multirow{5}{*}{2} & 10 & 48  & 0.020 & 0.003 & 0.149 & 0.0003 & 0.989\\ & 20 & 86 & 0.049 & 0.01 & 0.204 & 0.0121 & 0.983\\ & 50 & 196 & 0.191 & 0.033 & 0.172 & 0.0672 & 0.981\\  & 70 & 258 & 0.336 & 0.078 & 0.232 & 0.1332 & 0.989 \\ & 100 & 381 & 0.689 & 0.129 & 0.187 & 0.2412 & 0.995 \\
\hline

\multirow{5}{*}{3} & 10 & 1,438  & 0.750 & 0.159 & 0.212 & 0.0008 & 0.987\\ & 20 & 5,163 & 4.228 & 1.079 & 0.255 & 0.0243 & 0.986\\ & 50 & 32,558 & 53.96 & 12.09 & 0.224 & 0.1212 & 0.988\\  & 70 & 26,063 & 57.23 & 12.37 & 0.216 & 0.1872 & 0.985 \\ & 100 & 55,623 & 166.4 & 40.14 & 0.241 & 0.2799 & 0.994 \\
\hline

\multirow{5}{*}{4} & 10 & 88,812  & 60.36 & 11.20 & 0.185 & 0.012 & 0.967\\ & 20 & $1.05 \times 10^6$ & 1,226 & 91.73 & 0.074 & 0.066 & 0.994\\ & 50 & $4.39 \times 10^6$ & 10,258 & 2,886 & 0.281 & 0.2231 & 0.983\\ & 70 & $5.87 \times 10^6$ & 17,308 & 6,125 & 0.353 & 0.4215 & 0.984 \\ & 100 & $6.32 \times 10^6$ & 18,000+ & -- & -- & 0.749 & -- \\
\hline

\end{tabular}
\end{table}

Similar to the second set of experiments, it can be observed from Table \ref{tab:1D-PMCLP-PC-QoS} that $T_H$ (time taken by greedy approximation method for 1D-PMCLP-PC-QoS) is always less than $1$ second. The last column, $\alpha$, refers to the ratio of reward captured by the greedy approach to the optimal covered reward. Table \ref{tab:1D-PMCLP-PC-QoS} represents that $\alpha$ is at least $0.967$. Therefore, Table \ref{tab:1D-PMCLP-PC-QoS} represents the existing trade off between computational time and optimal/suboptimal coverage in the 1D-PMCLP-PC-QoS. It can help in deciding whether to use greedy or exact algorithm for a specific instances category depending on requested accuracy and available computational resources.

\section{Conclusion} \label{sec:Conclusion} \vspace{-0em}

We introduced a new generalization of the classical planar maximum coverage location problem where demand zones (DZs) and service zone (SZ) of facilities are represented by two-dimensional spatial objects (e.g., polygons, circles, etc.), the partial coverage is allowed in its true sense, facilities are allowed to be located anywhere on the continuous plane, and each facility has adjustable service range or quality of service (QoS). We denoted this problem by PMCLP-PC-QoS.  
We presented greedy and psuedo-greedy algorithms for the PMCLP-PC-QoS that have approximation ratio of $1 - 1/e$ and $1 - 1/e^{\eta}$ for $\eta\leq 1$, respectively. We investigated theoretical properties of the objective function of PMCLP-PC-QoS with rectangular DZs and SZs, denoted by PMCLP-PCR-QoS, and reduced the solution search space. We also proposed exact branch-and-bound based algorithm for the  PMCLP-PCR-QoS and one-dimensional PMCLP-PC-QoS (1D-PMCLP-PC-QoS) with SZs of different dimensions. 
For PMCLP-PCR-QoS with fixed and same QoS ($m=1$), the proposed algorithm is faster than the existing benchmark algorithm \citep{BanKia17IJOC}. Based on the computational results, we observed that the proposed greedy and exact algorithms for PMCLP-PCR-QoS and 1D-PMCLP-PC-QoS 
are computationally efficient and effective as well. 


%
%
%
%




{\bf Acknowledgments.} {This research is funded by Automotive Research Center (ARC) in accordance with Cooperative Agreement W56HZV-19-2-0001 U.S. Army CCDC Ground Vehicle Systems Center (GVSC) Warren, MI and by National Science Foundation Grant CMMI– 1824897, which are gratefully acknowledged.}




\bibliographystyle{informs2014} 
\bibliography{bibGenMCLP,bibGenMCLP_PartialCoverage,References} 

\end{document}